\documentclass[12pt]{article}
\nonstopmode
\RequirePackage[colorlinks,citecolor=blue,urlcolor=blue,linkcolor=blue]{hyperref}
\usepackage{graphicx,xspace,colortbl}
\usepackage{amsmath,amsthm,amsfonts}
\usepackage{color}
\usepackage{pdfsync}
\hypersetup{
colorlinks = true,
citecolor=blue,
urlcolor=blue,
linkcolor=blue,
pdfauthor = {A. Kuznetsov, J. C. Pardo, M. Savov},
pdfkeywords = {60G51,  Levy processes, exponential functional, integral equations, Mellin transform, asymptotic expansions},
pdftitle = {Distributional properties  of exponential functionals of Levy processes},
pdfpagemode = UseNone
}
\newtheorem{theorem}{Theorem}
\newtheorem{lemma}{Lemma}
\newtheorem{proposition}{Proposition}
\newtheorem{corollary}{Corollary}
\theoremstyle{definition}

\newtheorem{assumption}{Assumption}
\newtheorem{remark}{Remark}

\numberwithin{equation}{section}
\newcommand{\LL}{L\'{e}vy }

\newcommand{\PPPn}{\overline{\overline{\Pi}}^{(-)}}
\newcommand{\PPPp}{\overline{\overline{\Pi}}^{(+)}}
\newcommand{\tto}{_{t\to 0}}

\newcommand{\LLU}{{\mathfrak{L}}^{(U)}f(x)}
\newcommand{\IInf}{\int_{0}^{\infty}}
\newcommand{\lb}{\left (}
\newcommand{\rb}{\right )}
\newcommand{\lbb}{\left [}
\newcommand{\rbb}{\right ]}
    \def\beq{\begin{eqnarray}}
    \def\eeq{\end{eqnarray}}
    \def\beqq{\begin{eqnarray*}}
    \def\eeqq{\end{eqnarray*}}
    \def\c{{\mathbb C}}
    \def\p{{\mathbb P}}
    \def\e{{\mathbb E}}
    \def\re{\textnormal {Re}}
    \def\im{\textnormal {Im}}
    \def\i{\textnormal {i}}
    \def\d{{\textnormal d}}
    \def\r{{\mathbb R}}
    
    \def\mm{{\mathcal{M}}}
\title{
\textbf{Distributional properties  of exponential functionals of L\'evy processes}
}
\author{
\textbf{
A. Kuznetsov
\footnote{Department of Mathematics and Statistics,
York University,
4700 Keele Street,
Toronto, Ontario,
M3J 1P3, Canada. Email: {\color{blue} kuznetsov@mathstat.yorku.ca}. Research supported by the
Natural Sciences and Engineering Research Council of Canada. }
}
,
\,
\textbf{J. C. Pardo
\footnote{Centro de Investigaci\'on en Matem\'aticas A.C. Calle Jalisco s/n. 36240 Guanajuato, M\'exico. Email:
 {\color{blue} jcpardo@cimat.mx}. Research supported by CONACYT.}
}
,
\,
\textbf{
M. Savov
\footnote{New College,
Holywell Street,
Oxford, OX1 3BN,
UK. Email: {\color{blue} savov@stat.ox.ac.uk} ;  {\color{blue} mladensavov@hotmail.com}}
}
}
\date{\footnotesize This version: \today}

\begin{document}
\maketitle
\begin{abstract}
We study the distribution of the exponential functional $I(\xi,\eta)=\int_0^{\infty} \exp(\xi_{t-}) \d \eta_t$, where $\xi$ and $\eta$ are independent L\'evy processes. In the general setting, using the theory of Markov processes and Schwartz distributions, we prove that the law of this exponential functional satisfies an integral equation, which generalizes Proposition 2.1 in \cite{CPY}. In the special case when $\eta$ is a Brownian motion  with drift, we show that 
this integral equation leads to an important functional equation for the Mellin transform of $I(\xi,\eta)$, which 
  proves to be a very useful tool for studying the distributional properties of this random variable. For general L\'evy process
$\xi$ ($\eta$ being Brownian motion with drift) we prove that the exponential functional has a smooth density on $\r \setminus \{0\}$, but surprisingly the second derivative at zero may fail to exist. Under the additional assumption that $\xi$ has some positive
exponential moments we establish an asymptotic behaviour of $\p(I(\xi,\eta)>x)$ as $x\to +\infty$, 
and under similar assumptions on the negative exponential moments of $\xi$ we obtain a precise asymptotic expansion of the density of $I(\xi,\eta)$ as $x\to 0$. Under further assumptions on the L\'evy process $\xi$ one is able to prove much stronger results about the 
density of the exponential functional and we illustrate some of the ideas and techniques for the case when $\xi$ has hyper-exponential jumps.

\bigskip

\noindent {\it Keywords:}  L\'evy processes, exponential functional, integral equations, Mellin transform, asymptotic expansions.

\medskip

\noindent{\it AMS 2000 subject classifications: 60G51.}

\bigskip

\noindent Submitted to EJP on June 30 2011,  final version accepted December 1 2011.
\end{abstract}

\newpage

\section{Introduction}\label{section_introduction}

In this paper, we are  interested in  studying  distributional properties  of the random variable
\begin{equation}\label{I}
I(\xi,\eta):=\int_{0}^{\infty}e^{\xi_{t-}}\mathrm{d}\eta_{t},
\end{equation}
where $\xi$ and $\eta$ are independent real-valued \LL processes such that $\xi$ drifts to $-\infty$ and $\e[|\xi_{1}|]<\infty$ and $\e[|\eta_{1}|]<\infty$.
 
 The exponential functionals $I(\xi,\eta)$ appear in various aspects of probability theory. They describe the stationary measure of generalized Ornstein-Uhlenbeck processes and the entrance law of positive self-similar Markov processes, see \cite{BYS, CPY}. They also play a role in the theory of fragmentation processes and  branching processes, see \cite{Bfrag, KP}. Besides their theoretical value, the exponential functionals are very important objects in Mathematical Finance and Insurance Mathematics. They are related to Asian options, present values of certain perpetuities, etc., see  \cite{Duf, KLM, GP97} for some particular examples and results.
 
  In general, the distribution of exponential functionals is difficult to study. It is known explicitly only in some very special cases, see \cite{CaiKou2010, GP97, Ku11}. Properties of the distribution of $I(\xi,\eta)$ are also of particular interest. Lindner and Sato \cite{LS09} show that the density of $I(\xi,\eta)$ doesn't always exist, and in the special case when $\xi$ and $\eta$ are specific compound Poisson processes,
  distributional properties of $I(\xi,\eta)$ can be related to the problem of absolute continuity of the distribution of Bernoulli convolutions, which dates back to Erd\H os, see \cite{Erd40}. The distribution of $I(\xi,\eta)$, when $\xi_{s}=-s$ 
and in some other instances, is known to be self-decomposable and hence absolutely continuous,  see \cite{BLM, KMS06}. When $\eta$ is a subordinator with a strictly positive drift, the law of the exponential functional  $I(\xi, \eta)$ is absolutely continuous, see Theorem 3.9 in Bertoin et al. \cite{BLM}. Some further results are obtained in \cite{LA, PP, Pau, yorbook}.
  
The asymptotic behaviour $\p \lb I(\xi,\eta)> x \rb$, as $x\to\infty$, is a question which has attracted the attention of many researchers. In the general case, but under rather stringent requirements on the existence of exponential moments for $\xi$ and absolute moments for $\eta$, it has been studied in \cite{LindnerMaller2005}. The special case when $\eta_{t}=t$ has been considered in \cite{MZ, RI1, r2007} and properties of the density of the law of $I(\xi,\eta)$ at zero and infinity have been studied by \cite{ Ku11, KUPA11, PRS} and results such as asymptotic and convergent series expansions for the density have been obtained.

The first objective of this paper is to develop a general integral equation for the law of $I(\xi,\eta)$ under the assumptions that 
 $\e[|\xi_{1}|]<\infty$, $\e[\xi_1]<0$, $\e[|\eta_{1}|]<\infty$ and $\xi$ being independent of $\eta$. Using the fact that in general $I(\xi,\eta)$ is a stationary law of a generalized Ornstein-Uhlenbeck process, Carmona et al. \cite{CPY} show that if $\xi$ has jumps of bounded variation and $\eta_{t}=t$ then 
 the law of $I(\xi,\eta)$ satisfies a certain integral equation. We refine and strengthen their approach and using both stationarity properties of $I(\xi,\eta)$ and Schwartz theory of distributions, we show that in the general setting the law of $I(\xi,\eta)$ satisfies a certain integral equation. This equation is important on its own right, as demonstrated by Corollary \ref{Corollary10}, but it is also amenable to different useful transformations as can be seen from the discussion below.

The second main objective of the paper is to study some properties of $I_{\mu,\sigma}:=I(\xi,\eta)$ 
in the specific case when $\eta_{s}=\mu s+\sigma B_{s}$, where $B_{s}$ is a standard Brownian motion. 
Quantities of this type have already appeared in the literature, see \cite{GP97}, but have not been thoroughly studied. The latter, as it seems to us, is due to the lack of suitable techniques, which are available in the case when $\eta_{s}=s$, and in particular due to the lack of
any information about the Mellin transform of $I(\xi,\eta)$, which is the key tool for studying the properties of $I(\xi,\eta)$, see \cite{ Ku11, KUPA11, MZ}. We use the integral equation \eqref{SS1} and combine techniques from special functions, complex analysis and probability theory   to study the Mellin transform of $I_{\mu,\sigma}$, which is defined as $\mathcal{M}(s)=\e \lbb (I_{\mu,\sigma})^{s-1}\textbf{1}_{\{I_{\mu,\sigma}>0\}}\rbb$. In particular we derive an important functional equation for $\mathcal{M}(s)$, see \eqref{RecurBr}, and study the decay of $\mathcal{M}(s)$ 
as $\im(s) \to \infty$. These results supply us with quite powerful tools for studying the properties of the density of $I_{\mu,\sigma}$ via the Mellin inversion. Furthermore, the functional equation \eqref{RecurBr} allows for a meromorphic extension of $\mathcal{M}(s)$ when $\xi$ has some exponential moments. This culminates in very precise asymptotic results for $\p\lb I_{\mu,\sigma}>x \rb$, as $x\rightarrow\infty$, see Theorem \ref{thm_asym_infty}, and asymptotic expansions for $k(x)$, the density of $I_{\mu,\sigma}$, as $x\to 0$, see Theorem \ref{Asym0}. The latter results show us that while $k(x)\in C^{\infty}(\r\setminus\{0\})$, rather unexpectedly $k''(0)$ may not exist. Finally, we would like to point out that while the behaviour of $\p\lb I_{\mu,\sigma}>x \rb$, as $x\rightarrow\infty$, might be partially studied via the fact that $I_{\mu,\sigma}$ solves a random recurrence equation, see for example \cite{LindnerMaller2005}, the behaviour of $k(x)$, as $x\to 0$, seems for the moment to be
only tractable via our approach based on the Mellin transform.

As another illustration of possible applications of our general results, we study the density  of $I_{\mu,\sigma}$  when $\xi$ has hyper-exponential jumps (see \cite{Cai2009127, CaiKou2010, KUKYPA10}). This class of processes is quite important for applications 
in Mathematical Finance and Insurance Mathematics, and it is particularly well suited for investigation using our methods due to the rich 
analytical structure enjoyed by these processes. In this case we show how to derive complete asymptotic expansions 
of $k(x)$ both at zero and infinity. We point out that our methodology is not restricted to this particular case, and can be easily applied to more general classes of L\'evy processes.

The paper is organized as follows: in Section \ref{section_integral_equations}, we study the law of $I(\xi,\eta)$
 for general independent \LL processes $\xi$ and $\eta$ and derive an integral equation for the law of $I(\xi,\eta)$; in Section \ref{section_BM_exp_functionals}, we specialize the results obtained in Section \ref{section_integral_equations} to the case when $\eta_{s}=\mu s+\sigma B_{s}$ and, employing additionally various techniques from special functions and complex analysis, we study the properties of the density of $I_{\mu,\sigma}$. Section \ref{section_applications} is devoted to some applications of the results derived in the previous section. In particular,  we study the asymptotic behaviour at infinity of the tail of $I_{\mu,\sigma}$ and of its density at zero, and in the case of processes with hyper-exponential jumps, we show how these results can be considerably strengthened.

\section{Integral equation satisfied by the law of $I(\xi,\eta)$}\label{section_integral_equations}

Let us introduce some notation which will be used throughout this paper. The main underlying objects are two independent L\'evy processes
$\xi$ and $\eta$ defined on a probability space $(\Omega,{\mathcal F},\p)$. As is standard, we assume that both processes are started from zero under the probability  measure $\p$. 

\begin{assumption}\label{assumption1}
Everywhere in this paper we will assume that  
\beq\label{assumption_first_moments}
\e[|\xi_{1}|]<\infty, \;\;\; \e[\xi_{1}]<0, \;\;\;  \e[|\eta_{1}|]<\infty.
\eeq 
\end{assumption}

The  characteristics of the L\'evy processes $\xi$ and $\eta$ will be denoted by $(b_{\xi},\sigma_{\xi}, \Pi_{\xi})$ and 
$(b_{\eta},\sigma_{\eta}, \Pi_{\eta})$. In particular $\Pi_\xi(\d x)$ and $\Pi_\eta(\d x)$ are the \LL measures of $\xi$ and $\eta$, respectively. 
We use the following notation for the double-integrated tail
\[
\PPPp_\xi(x)=\int_{x}^{\infty}\Pi_\xi((y,\infty))\d y\qquad\text{ and }\qquad\PPPn_\xi(x)=\int_{x}^{\infty}\Pi_\xi((-\infty,-y))\d y,
\]
and similarly for  $\PPPp_\eta$ and $\PPPn_\eta$. Using the L\'evy-It\^o decomposition (see Theorem 2.1 in \cite{Kyprianou}) it is easy to check that Assumption \ref{assumption1}  implies that the above quantities are finite for all $x>0$.

We define the Laplace exponents $\psi_{\xi}(z)=\ln{\lb\e\lbb e^{z\xi_{1}}\rbb\rb}$ and $\psi_{\eta}(z)=\ln{\lb\e\lbb e^{z\eta_{1}}\rbb\rb}$, where without any further assumptions $\psi_{\xi}$ and $\psi_{\eta}$ are defined at least for $\re(z)=0$, see \cite[Chapter I]{Be}. 
The Laplace exponent $\psi_{\xi}$ can be expressed in the following two equivalent ways
 \beq\label{SpecialFormExponents}
 \psi_{\xi}(z)&=&\frac{\sigma^{2}_{\xi}}{2}z^{2}+b_{\xi}z+\int_{\r} \left( e^{zx}-1-zx\right) \Pi_{\xi}(\d x)\\ \nonumber
&=&\frac{\sigma^{2}_{\xi}}{2}z^{2}+b_{\xi}z+z^{2}\lb\IInf\PPPp_{\xi}(w)e^{xz}dx+\IInf\PPPn_{\xi}(x)e^{-xz}dx\rb,
 \eeq
with a similar expression for $\psi_{\eta}$. The first equality in \eqref{SpecialFormExponents} is essentially 
the L\'evy-Khintchine  formula (see Theorem 1 in \cite{Be}) with the cutoff function $h(x)\equiv 1$. The standard choice for the cutoff function
in the L\'evy-Khintchine formula would be ${\mathbf 1}_{\{|x|<1\}}$, however it is well-known that if $\e[|\xi_1|]<\infty$ then we can take a simpler cutoff function $h(x)\equiv 1$. The second equality in \eqref{SpecialFormExponents} follows easily by repeated integration by parts. 
Note that according to \eqref{SpecialFormExponents}, we have $b_{\xi}=\psi_{\xi}'(0)=\e[\xi_1]$ and similarly 
$b_{\eta}=\psi_{\eta}'(0)=\e[\eta_1]$. 

We recall that the exponential functional $I(\xi,\eta)$ is defined by \eqref{I}, its law will be denoted by $m(\d x):=
\p(I(\xi,\eta)\in \d x)$. The density of $I(\xi,\eta)$, provided it exists, will be denoted by $k(x)$.

 Our main result in this section is the derivation of an integral equation for the law of $I(\xi,\eta)$. 
 This equation will be very useful  later, when we'll derive the functional equation \eqref{RecurBr}
 for the Mellin transform of the exponential functional in the special case when $\eta$ is a Brownian motion with drift.   
 The main idea of this Theorem comes from Proposition 2.1 in \cite{CPY}.
 
\begin{theorem}\label{Theorem_main_integral_equation}
Assume that condition \eqref{assumption_first_moments} is satisfied. Then the exponential functional $I(\xi,\eta)$ is well defined and its law satisfies the following integral equation: for $v>0$
\beq\label{SS1}
\nonumber
&&\lb b_{\xi}\int_{v}^{\infty}m(\d x)\rb\d v+\frac{\sigma_{\xi}^{2}}{2}vm(\d v)+\lb\int_{v}^{\infty}\PPPn_\xi\left(\ln{\frac{x}{v}}\right)m(\d x)\rb\d v+\lb\int_{0}^{v}\PPPp_\xi\left(\ln{\frac{v}{x}}\right)m(\d x)\rb\d v\\ \nonumber
&&+\lb b_{\eta}\int_{v}^{\infty}\frac{m(\d x)}{x}\rb\d v+\frac{\sigma_{\eta}^{2}}{2}\frac{m(\d v)}{v}-\lb\frac{\sigma_{\eta}^{2}}{2}\int_{v}^{\infty}\frac{m(\d x)}{x^{2}}\rb\d v\\
&&+\lb\frac{1}{v}\int_{0}^{v}\PPPp_{\eta}(v-x)m(\d x)\rb\d v+\lb\frac{1}{v}\int_{v}^{\infty}\PPPn_{\eta}(x-v)m(\d x)\rb\d v\\ \nonumber
&&-\lb\int_{v}^{\infty}\frac{1}{w^{2}}\int_{0}^{w}\PPPp_{\eta}(w-x)m(\d x)\d w\rb\d v-\lb\int_{v}^{\infty}\frac{1}{w^{2}}\int_{w}^{\infty}\PPPn_{\eta}(x-w)m(\d x)\d w\rb\d v=0,
\eeq
 where all quantities in \eqref{SS1} are a.e. finite. Equation \eqref{SS1} for the law of $I(\xi,-\eta)$ on $(0,\infty)$  describes $m(\d x)$ on $(-\infty,0)$.  
\end{theorem}

\vspace{0.25cm}
{\color{red} {\bf Remark (added in June 2020).}  As was shown by Behme et. al. in \cite{BeLiRe}[Remark 5.10],  the equation 
\eqref{SS1} is valid only under an additional assumption: $\eta$ must be (i) a subordinator or (ii) a process without positive jumps. We would like to thank the authors of \cite{BeLiRe} for pointing this out. Note, that the equation \eqref{SS1} is used in Sections \ref{section_BM_exp_functionals} and \ref{section_applications} 
only in the latter case:  
$\eta$ is a Brownian motion with drift and thus it has no positive jumps, thus all of our results in Sections \ref{section_BM_exp_functionals} and \ref{section_applications} remain valid.}
\vspace{0.25cm}

The proof of Theorem \ref{Theorem_main_integral_equation} is based on the 
so-called generalized Ornstein-Uhlenbeck (GOU) process, which is defined as
\begin{equation}\label{Orn}
U_{t}=U_{t}(\xi,\eta)=xe^{\xi_{t}}+e^{\xi_{t}}\int_{0}^{t}e^{-\xi_{s-}}\d\eta_{s}\stackrel{d}=xe^{\xi_{t}}+\int_{0}^{t}e^{\xi_{s-}}\d\eta_{s},\text{ for  $t>0$.}
\end{equation}
 Note that the GOU process is a strong Markov process, see \cite[Appendix 1]{CPY}. Lindner and Maller \cite{LindnerMaller2005} have shown that the existence of a stationary distribution for the GOU process is closely related to the a.s. convergence of the stochastic integral $\int_{0}^{t}e^{\xi_{s-}}\d\eta_{s}$, as $t\to \infty$. Necessary and  sufficient conditions for the convergence of $I(\xi,\eta)$ were obtained by Erickson and Maller \cite{EM05}. More precisely, they showed that this happens if and only if
\begin{equation}\label{condEM}
\lim_{t\to\infty}\xi_t=-\infty\;\;\;\textrm{ a.s.} \qquad\textrm{ and } \qquad \int_{\r\setminus[-e,e]}\left[\frac{\log|y|}{1+\int_1^{\log|y|\vee 1}\Pi_{\xi}(\r \setminus (-z,z))\d z}\right]\Pi_\eta(\d y)<\infty.
\end{equation}
It is easy to see that Assumption \ref{assumption1} implies (\ref{condEM}). Hence $I(\xi,\eta)$ is well-defined and the stationary distribution satisfies  $U_\infty\stackrel{d}=I(\xi,\eta)$. This identity in distribution is the starting point of the proof of Theorem \ref{Theorem_main_integral_equation}.

As the proof of Theorem \ref{Theorem_main_integral_equation} is rather long and technical, we will divide it into several steps. We first compute the generator of $U$, here denoted by $\mathfrak{L}^{(U)}$. This result may be of independent interest, therefore we present it in Proposition \ref{Prop1} below. Then we note that the stationary measure $m(\d x)$ satisfies the equation
\begin{equation}\label{Schwarz}
\int_{0}^{\infty}\mathfrak{L}^{(U)}f(x)m(\d x)=0,
\end{equation}
 where $f$ is any infinitely differentiable function with a compact support in $(0,\infty)$. Indeed, \eqref{Schwarz}  follows from (2.1) in \cite{CPY} or from the definition of infinitesimal generator and the observation that, for all $t\geq0$,
\[\int_{0}^{\infty}\e[f(U_t)]m(dx)=\int_{0}^{\infty}f(x)m(dx).\] 
  Finally, an application of Schwartz theory of distributions after rephrasing \eqref{Schwarz} gives \eqref{SS1}.  

We start by working out how the infinitesimal generator  of $U$, i.e. $\mathfrak{L}^{(U)}$, acts on functions in $\mathcal{K}\subset C_{0}(\r)$, where
\beq\label{core} \nonumber
\mathcal{K}&=&\left\{f(x):\,f(x)\in C^{2}_{b}(\mathbb{R}),\,f(e^{x})\in C^{2}_{b}(\mathbb{R})\cap C_{0}(\r)\right\}
\\ \qquad \qquad\qquad && 
\cap\left\{\,f(x)=0,\,\text{ for $x\leq0$; }f'(0)=f''(0)=0\right\}
\eeq
and $C^{2}_{b}(\r)$ stands for two times differentiable, bounded functions with bounded derivatives on $\r$ and $C_{0}(\r)$ is the set of continuous functions vanishing at $\pm\infty$. Denote by $\mathfrak{L}^{(\xi)}$ and $\mathfrak{L}^{(\eta)}$ ( resp. $\mathfrak{D}^\xi$ and $\mathfrak{D}^\eta$) the infinitesimal generators (resp. domains) of $\xi$ and $\eta$. Note that
\beq\label{MoreGeneretors}
\mathfrak{L}^{(\xi)}f(x)&=&b_{\xi}f'(x)+\frac{\sigma^{2}_{\xi}}{2}f''(x)+\int_{\r}\lb f(x+y)-f(x)-yf'(x)\rb\Pi_{\xi}(dy)\\ \nonumber
 &=&b_{\xi}f'(x)+\frac{\sigma^{2}_{\xi}}{2}f''(x)+\int_{\r_+}f''(x+w)\PPPp_{\xi}(w)dw+\int_{\r_+}f''(x-w)\PPPn_{\xi}(w)dw,
\eeq
with a similar expression for $\mathfrak{L}^{(\eta)}$. The first formula in \eqref{MoreGeneretors} is a trivial modification of the form of the generator of \LL processes for the case when the cutoff function is $h(x)\equiv 1$, see \cite[p. 24]{Be}, whereas the second expression follows easily by integration by parts, the fact that $f\in\mathcal{K}$ and $\e [|\xi_{1}|]<\infty$.
Finally, we are ready to state our result, which should strictly be seen as an extension of Proposition 5.8 in \cite{CPY} where the generator $\mathfrak{L}^{(U)}$ has been derived under very stringent conditions.

\begin{proposition}\label{Prop1}
Assume that condition \eqref{assumption_first_moments} is satisfied. Let $f\in\mathcal{K}$, $g(x):=(xf'(x))$ and $\phi(x):=f(e^{x})$. Then, $f\in \mathfrak{D}^\eta$, $\phi\in \mathfrak{D}^{\xi}$ and
\begin{eqnarray}\label{generator}
\nonumber&&\mathfrak{L}^{(U)}f(x)=\mathfrak{L}^{(\xi)}\phi(\ln{x})+\mathfrak{L}^{(\eta)}f(x)\\
\nonumber&&=b_{\xi}g(x)+\frac{\sigma_{\xi}^{2}}{2}xg'(x)+\int_{0}^{x}g'(v)\PPPn_{\xi}\Big(\ln{\frac{x}{v}}\Big)\d v+\int_{x}^{\infty}g'(v)\PPPp_{\xi}\Big(\ln{\frac{v}{x}}\Big)\d v\\
&&+b_{\eta}f'(x)+\frac{\sigma_{\eta}^{2}}{2}f''(x)+\IInf f''(x+w)\PPPp_{\eta}(w)\d w+\IInf f''(x-w)\PPPn_{\eta}(w)\d w.
\end{eqnarray}
\end{proposition}
\begin{proof}
The main idea is to use the definition of the infinitesimal generator and It\^o's formula. Let $f\in \mathcal{K}$ and note that by definition
\beqq
\LLU=\lim\tto\frac{\e_{x}\left[f(U_{t})\right]-f(x)}{t}=\lim\tto\frac{1}{t}\left(\e\left[f\left(xe^{\xi_{t}}+\int_{0}^{t}e^{\xi_{s-}}\d\eta_{s}\right)\right]-f(x)\right).
\eeqq
Using the fact that $\lb U_t\rb_{t\geq 0}$ is a semimartingale and $f\in\mathcal{K}$, we apply It\^o's formula to $f(U_t)$  to obtain
\begin{equation}\label{itof}
f(U_{t})-f(x)=\int_{0}^{t}f'(U_{s-})\d U_{s}+\frac{1}{2}\int_{0}^{t}f''(U_{s-})\d [U,U]^{c}_{s}+\sum_{s\leq t}\lb f(U_{s})-f(U_{s-})-\Delta U_{s}f'(U_{s-})\rb.
\end{equation}
Now, let
$H_t:=e^{\xi_{t}}$ and $V_t:=x+\int_{0}^{t}e^{-\xi_{s-}}d\eta_{s},$
and note that $U_t=H_tV_t$. Hence by integration by parts
\begin{align*}
U_{t}=x+\int_{0}^{t}H_{s-}\d V_{t}+\int_{0}^{t}V_{s-}\d H_{s}+[H,V]_t.
\end{align*}

Using the L\'evy-It\^o decomposition (see Theorem 2.1 in \cite{Kyprianou}) and Assumption \ref{assumption1}, we find that the L\'evy processes $\xi$ and $\eta$ can be written as follows
\begin{equation}\label{Levy}
\xi_{t}=\sigma_\xi B_{t}+b_\xi t+X_{t}, \qquad \eta_t=\sigma_\eta W_{t}+b_\eta t+Y_{t},
\end{equation}
 where $B$ and $W$ are Brownian motions, $X$ and $Y$ are pure jump zero mean martingales, and the  processes $B, W, X$ and $Y$ are mutually independent.  Then we get
\[
V_{t}=x+b_{\eta}\int_{0}^{t}e^{-\xi_{s-}}\d s+\sigma_{\eta}\int_{0}^{t}e^{-\xi_{s-}}\d W_{s}+N_{t},
\]
where $N_{t}=\int_{0}^{t}e^{-\xi_{s-}}\d Y_{s}$ is a pure jump local martingale. On the other hand using It\^o's formula, we have
\[
\begin{split}
H_{t}&=e^{\xi_t}=1+\int_{0}^{t}e^{\xi_{s-}}\d\xi_{s}+\frac{1}{2}\int_{0}^{t}e^{\xi_{s-}}\d[\xi,\xi]^{c}_{s}+\sum_{s\leq t}e^{\xi_{s-}}(e^{\Delta \xi_{s}}-\Delta \xi_{s}-1)\\
&=1+\left(b_\xi+\frac{\sigma^{2}_{\xi}}{2}\right)\int_{0}^{t}e^{\xi_{s-}}\d s+\sigma_{\xi}\int_{0}^{t}e^{\xi_{s-}}\d B_{s}+\widetilde{N}_{t}+\sum_{s\leq t}e^{\xi_{s-}}\Big(e^{\Delta \xi_{s}}-\Delta \xi_{s}-1\Big),
\end{split}
\]
where $\widetilde{N}_{s}=\int_{0}^{t}e^{\xi_{s}}\d X_{s}$ is a pure jump local martingale. Therefore, we conclude that
\begin{align*}
[H,V]_t=\left[\sigma_{\xi}\int_{0}^{t}e^{-\xi_{s-}}\d B_{s},\sigma_{\eta}\int_{0}^{t} e^{-\xi_{s-}}\d W_{s}\right]_t +\sum_{s\le t} \Delta V_s\Delta H_s=0\qquad \textrm{a.s.,}
\end{align*}
since $\Delta V_s=e^{-\xi_{s-}} \Delta \eta_s$, $\Delta H_s=H_{s-}(e^{\Delta \xi_s}-1)$ and the fact that $\xi$ and $\eta$ are independent and do not jump simultaneously a.s.  This implies that 
\[
U_{t}=x+\int_{0}^{t}H_{s-}\d V_{s}+\int_{0}^{t}V_{s-}\d H_{s}=x+\int_{0}^{t}e^{\xi_{s-}}\d V_{s}+\int_{0}^{t}V_{s-}\d H_{s}.
\]
Using the expressions of $H$ and $V$, we deduce that
\[
\begin{split}
U_{t}&=x+b_{\eta}t +\sigma_{\eta}W_t+\int_0^t e^{\xi_{s-}}\d N_s+\left(b_{\xi}+\frac{\sigma^{2}_{\xi}}{2}\right)\int_{0}^{t}V_{s-}e^{\xi_{s-}}\d s\\
&\hspace{2cm}+\sigma_{\xi}\int_{0}^{t}V_{s-}e^{\xi_{s-}}\d B_{s}+\int_0^t V_{s-}\d \widetilde{N}_s+\sum_{s\leq t}V_{s-}e^{\xi_{s-}}\Big(e^{\Delta \xi_{s}}-\Delta \xi_{s}-1\Big)\\
&=x+K_{t}+K^{c}_{t}+b_{\eta}t+\left(b_{\xi}+\frac{\sigma^{2}_{\xi}}{2}\right)\int_{0}^{t}U_{s-}\d s+\sum_{s\leq t}U_{s-}\lb e^{\Delta \xi_{s}}-\Delta \xi_{s}-1\rb,
\end{split}
\]
where
\beqq
K_t&=&\int_0^t e^{\xi_{s-}}\d N_s+\int_0^t V_{s-}\d \widetilde{N}_s
=Y_t+\int_0^t U_{s-} d X_s, \\
 K^{c}_{t}&=&\sigma_{\eta}W_t+\sigma_{\xi}\int_{0}^{t}U_{s-}\d B_{s}.
\eeqq
From the  definition  of $K$ and $K^c$, and the mutual independence of $B$, $W$, $N$ and $\widetilde{N}$, we get for the continuous part of the quadratic variation of $U$
\begin{align*}
[U,U]^{c}_t=[K^{c},K^{c}]_t=\sigma^{2}_{\eta} t+\sigma^{2}_{\xi}\int_0^t U^2_{s-} \d s.
\end{align*}
 Putting all the pieces  together in  identity (\ref{itof}), we have
\[
\begin{split}
f(U_{t})-f(x)&=M_t+b_{\eta}\int_{0}^{t}f'(U_{s-})\d s+\left(b_{\xi}+\frac{\sigma^{2}_{\xi}}{2}\right)\int_{0}^{t}f'(U_{s-})U_{s-}\d s\\
&+\sum_{s\leq t}f'(U_{s-})U_{s-}\Big(e^{\Delta \xi_{s}}-\Delta \xi_{s}-1\Big)+\frac{\sigma^{2}_{\eta}}{2}\int_{0}^{t}f''(U_{s-})\d s+\frac{\sigma^{2}_{\xi}}{2}\int_{0}^{t}f''(U_{s-})U^{2}_{s-}\d s\\
&+\sum_{s\leq t}\left(f(U_{s})-f(U_{s-})-\Delta U_{s}f'(U_{s-})\right)
\end{split}
\]
where $M$ is a local martingale starting from $0$ and $M$ describes the integration with respect to $K$ and $K^{c}$ in the expressions above. Using the fact that $f\in\mathcal{K}$ implies $f(x)=0$ for $x<0$ and $x|f'(x)|+x^{2}|f''(x)|<C(f)<\infty$, we deduce that $M_{t}$ is a proper martingale as all other terms in the expression above have a finite absolute first moment. Furthermore applying the compensation formula to the jump part of $f(U_t)$ we get
\[
\e\left[\sum_{s\leq t}f'(U_{s-})U_{s-}\lb e^{\Delta \xi_{s}}-\Delta \xi_{s}-1\rb\right]=\e\left[\int_{0}^{t}f'(U_{s-})U_{s-}\left(\int_{ y\in \r}(e^{y}-y-1)\Pi_{\xi}(dy)\right) \d s\right].
\]
Similarly, using the fact that $\Delta U_{s}=\Delta \eta_{s}$ when $\Delta\eta_{s}\neq0$ and $\Delta U_{s}=U_{s-}(e^{\Delta \xi_{s}}-1)$ when $\Delta\xi_{s}\neq0$ (see the definition of $U$) we get
\beqq
\e\lbb\sum_{s\leq t}\lb f(U_{s})-f(U_{s-})-\Delta U_{s-}f'(U_{s-})\rb\rbb  =\e\left[\int_{0}^{t}\int_{z\in\mathbb{R}}\left(f(U_{s-}+z)-f(U_{s-})-zf'(U_{s-})\right)\Pi_{\eta}(\d z)\d s\right]\\
+\e\left[\int_{0}^{t}\int_{y\in\r}\lb f(U_{s-}e^{y})-f(U_{s-})-\Big(e^{y}-1\Big)f'(U_{s-})U_{s-}\rb\Pi_{\xi}(\d y)\d s\right].
\eeqq
Finally, as $f\in\mathcal{K}$, we derive
\[
\begin{split}
\e\Big[f(U_{t})\Big]-f(&x)=b_{\eta}\e\left[\int_{0}^{t}f'(U_{s-})\d s\right]+\left(b_{\xi}+\frac{\sigma^{2}_{\xi}}{2}\right)\e\left[\int_{0}^{t}f'(U_{s-})U_{s-}\d s\right]+\frac{\sigma^{2}_{\eta}}{2}\e\left[\int_{0}^{t}f''(U_{s-})\d s\right]\\
&+\frac{\sigma^{2}_{\xi}}{2}\e\left[\int_{0}^{t}f''(U_{s-})U^{2}_{s-}\d s\right]+\e\left[\int_{0}^{t}\int_{z\in\mathbb{R}}\Big(f(U_{s-}+z)-f(U_{s-})-zf'(U_{s-})\Big)\Pi_{\eta}(\d z)\d s\right]\\
&+\e\left[\int_{0}^{t}\int_{y\in\r}\Big(f(U_{s-}e^{y})-f(U_{s-})-yf'(U_{s-})U_{s-}\Big)\Pi_{\xi}(\d y)\d s\right].
\end{split}
\]
and dividing by $t$, letting $t$ go to $0$ and recalling that $\tilde{U}_0=x$ a.s., we obtain for $f\in\mathcal{K}$ the identity
\begin{eqnarray}
\nonumber&&\mathfrak{L}^{(U)}f(x)=b_{\eta}f^\prime(x)+\left(b_{\xi}+\frac{\sigma^{2}_{\xi}}{2}\right)xf^\prime(x)+\frac{\sigma^{2}_{\eta}}{2}f^{\prime\prime}(x)+\frac{\sigma^{2}_{\xi}}{2}f^{\prime\prime}(x)x^2\\
&&+\int_{z\in\mathbb{R}}\Big(f(x+z)-f(x)-zf^\prime(x)\Big)\Pi_{\eta}(\d z)+\int_{y\in\r}\Big(f(xe^{y})-f(x)-yxf^\prime(x)\Big)\Pi_{\xi}(\d y),
\end{eqnarray}
 and  therefore the infinitesimal generator of $U$ satisfies
\begin{align*}
&\mathfrak{L}^{(U)}f(x)=\mathfrak{L}^{(\xi)}\phi(\ln{x})+\mathfrak{L}^{(\eta)}f(x).
\end{align*}
In order to finish the proof one only has to apply 
integration by parts. 
\end{proof}

\noindent The following Lemma will also be needed for our proof of Theorem \ref{Theorem_main_integral_equation}.
\begin{lemma}\label{auxLemma}
Assume that condition \eqref{assumption_first_moments} is satisfied. Let
$\nu(\d v)$ denote the measure in the left-hand side of formula \eqref{SS1}.
Then $|\nu|(\d v)$ and hence $\nu(\d v)$ define finite measures on any compact subset of $(0,\infty)$ and for any $a>0$
\begin{equation}\label{limitNU}
\lim_{z\to\infty}z^{-1}|\nu|\lb (a,z)\rb=0.
\end{equation}
\end{lemma}
\begin{proof}
We only need to prove \eqref{limitNU}, as the finiteness of $|\nu|(\d v)$ on compact subsets of $(0,\infty)$ follows from \eqref{limitNU}.
 It is sufficient to show the claims for $1\geq a>0$. We integrate every term on the left-hand side of \eqref{SS1} from $a$ to $z$ and then divide by $z$. This shows that the limit goes to zero, as $z\to\infty$. We first note that
\[
\lim_{z\to\infty}z^{-1}\int_{a}^{z}xm(\d x)=0\quad \textrm{ and }\quad \lim_{z\to\infty}z^{-1}\int_{a}^{z}\frac{m(\d x)}{x}\leq \lim_{z\to\infty}(az)^{-1}\int_{a}^{\infty}m(\d x)=0.
\]
Hence,
\[
\lim_{z\to \infty}z^{-1}\int_{a}^{z}\int_{v}^{\infty}m(\d x)\d v\le \lim_{z\to\infty}\left(z^{-1}\int_{a}^{z}x m(\d x)+\int_{z}^{\infty} m(\d x)\right)=0,
\]
\[
\lim_{z\to\infty}z^{-1}\int_{a}^{z}\int_{v}^{\infty}\frac{m(\d x)}{x}\d v=0\quad \textrm{ and }\quad
\lim_{z\to\infty}z^{-1}\int_{a}^{z}\int_{v}^{\infty}\frac{m(\d x)}{x^{2}}\d v=0.
\]
So far, we have checked that the terms in \eqref{SS1} that do not depend on the tail of the \LL measure vanish under the transformation we made, as $z\rightarrow \infty$. Now, we turn our attention to the terms that involve the \LL measure of $\xi$. When we'll be dealing with these integrals,
the main trick that we will use is to change the order of integration. First, we check that
\beqq
&&\limsup_{z\to\infty}z^{-1}\int_a^z\int_v^\infty\PPPn_{\xi}\left(\ln{\frac{x}{v}}\right)m(\d x)\d v \\
&&\le \limsup_{z\to\infty}z^{-1}\lb \int_a^z\int_v^{ev}\PPPn_{\xi}\left(\ln{\frac{x}{v}}\right)m(\d x)\d v\rb+\limsup_{z\to\infty}z^{-1}\lb \PPPn_{\xi}(1)\int_a^z m(ev,\infty)\d v\rb\\
&&=\limsup_{z\to\infty}z^{-1}\int_a^z\int_v^{ev}\PPPn_{\xi}\left(\ln{\frac{x}{v}}\right)m(\d x)\d v
\le \limsup_{z\to\infty}z^{-1}\int_a^{ez}\int_{x/e}^{x}\PPPn_\xi\left(\ln{\frac{x}{v}}\right)\d v\,m(\d x)\\
&&=  \left[\int_{0}^{1}\PPPn_{\xi}(w)e^{-w}\d w\right] \times \limsup_{z\to\infty}z^{-1}\int_a^{ez}x m(\d x)=0
\eeqq
where we have applied Fubini's Theorem, a change of variables $w=\ln(x/v)$ and we have used the finiteness of $\e[|\xi_{1}|]$ and henceforth the finiteness of the quantities  $\int_{0}^{1}\PPPn_\xi\big(w\big)\exp(-w)\d w$ and $\PPPp_{\xi}(1)$.
 
Next using Fubini's Theorem and the monotonicity of $\PPPp_\xi$, we note that for any positive number $b$,
\beqq
&&\limsup_{z\to\infty}z^{-1}\int_a^z\int_0^v\PPPp_\xi\bigg(\ln{\frac{v}{x}}\bigg)m(\d x)\d v\\ 
&&\le \limsup_{z\to\infty}z^{-1}\int_0^z\int_0^v\PPPp_\xi\left(\ln{\frac{v}{x}}\right)m(\d x)\d v
=\limsup_{z\to\infty} z^{-1}\int_0^z x\int_0^{\ln(z/x)}\PPPp_\xi\big(w\big)e^w\d w\,m(\d x)\\
&& \le \limsup_{z\to\infty}z^{-1}\lb \int_0^{b}\PPPp_\xi\big(w\big)e^w\d w\int_0^z x\,m(\d x)+ \int_0^z x\int_b^{\ln(z/x)\lor b}\PPPp_\xi\big(w\big)e^w\d w\,m(\d x)\rb\le \PPPp_\xi\big(b\big).
\eeqq
Since  $\PPPp(b)$  decreases to zero as $b$ increases, we see that
\[
\lim_{z\to \infty} z^{-1}\int_a^z\int_0^v\PPPp_{\xi}\left(\ln{\frac{v}{x}}\right)m(\d x)\,\d v=0.
\]
Since $\eta$ has a finite mean and $m$ is a finite measure
\beqq
&&\limsup_{z\to\infty}z^{-1}\int_a^z\frac{1}{v}\int_0^v\PPPp_\eta(v-x)m(\d x)\, \d v\\
&& \le \limsup_{z\to\infty}z^{-1}\lb \PPPp_\eta(a)\ln\lb\frac{z}{a}\rb+\int_a^z\frac{1}{v}\int_{v-a}^v\PPPp_\eta(v-x)m(\d x)\, \d v\rb \\
&&= \limsup_{z\to\infty}z^{-1}\lb \int_0^z m (\d x) \int_{a\lor x}^{(x+a )\land z} \PPPp_\eta(v-x)\, \frac{\d v}{v}\rb \\
&& \le \left[ \int_0^a \PPPp_\eta(s)\d s \right] \times
\lim_{z\to\infty}(az)^{-1} \int_0^z m (\d x) =0.
\eeqq
 Similarly, we estimate the following integral
\beqq
&&\limsup_{z\to\infty}z^{-1}\int_a^z\int_v^\infty \frac{1}{w^2}\int_0^w\PPPp_\eta(w-x)m(\d x)\, \d w\,\d v\\
&&\le\limsup_{z\to\infty}z^{-1}\lb \PPPp_\eta(a)\ln\lb\frac{z}{a}\rb+\int_a^z\int_v^\infty \frac{1}{w^2}\int_{w-a}^w\PPPp_\eta(w-x)m(\d x)\, \d w\,\d v\rb\\
&&= \limsup_{z\to\infty}z^{-1}\int_a^z\int_{v-a}^\infty \int_{v\lor x}^{x+a}\frac{1}{w^2}\PPPp_\eta(w-x) \d w\, m(\d x)\,\d v \\
&&\le\left[\int_0^a \PPPp_\eta(s)\d s\right] \times \limsup_{z\to\infty}z^{-1}\int_a^z\frac{1}{v^2}m(v-a,\infty)\,\d v=0.\\
\eeqq

As for the remaining two integrals, we split the innermost integrals at the point $x=v+a$ so that $\PPPn_\eta(x-v)=\PPPn_\eta(a)$ and similarly estimate the resulting two terms to get
\[\limsup_{z\to\infty}z^{-1}\int_a^z\frac{1}{v}\int_v^\infty\PPPn_\eta(x-v)m(\d x)\, \d v=\limsup_{z\to\infty}z^{-1}\int_a^z\int_v^\infty \frac{1}{w^2}\int_w^\infty\PPPn_\eta(x-w)m(\d x) \d w\d v=0.\]
Thus, we verify \eqref{limitNU} and conclude the proof of Lemma \ref{auxLemma}.
\end{proof}

\noindent Now that we have established Proposition \ref{Prop1} and Lemma \ref{auxLemma}, we are ready to complete the proof of Theorem  \ref{Theorem_main_integral_equation}.

\vspace{0.2cm}

\noindent{\it Proof of Theorem \ref{Theorem_main_integral_equation}}. 
Take an infinitely differentiable function $f$ with compact support in $(0,\infty)$ and let $g(x):=xf'(x)$. We use \eqref{Schwarz},  \eqref{generator}, and the identity  $g(x)=\int_{0}^{x}g'(v)\d v$ to get,
\[
\begin{split}
\IInf \mathfrak{L}^{(\xi)}\phi(\ln{x}) m(\d x)&=b_{\xi}\IInf g(x)m(\d x)+\frac{\sigma_{\xi}^{2}}{2}\IInf x g'(x)m(\d x)\\
&+\IInf \int_{0}^{x}g'(v)\PPPn_{\xi}\left(\ln{\frac{x}{v}}\right)\d vm(\d x)+\IInf\int_{x}^{\infty}g'(v)\PPPp_{\xi}\left(\ln{\frac{v}{x}}\right)\d vm(\d x)\\
&=\IInf g'(v)\lb b_{\xi}\int_{v}^{\infty}m(\d x)\rb\d v+\IInf g'(v)\lb\frac{\sigma_{\xi}^{2}}{2}v m(\d v)\rb \\
&+\IInf g'(v)\lb\int_{v}^{\infty}\PPPn_{\xi}\left(\ln{\frac{x}{v}}\right)m(\d x)\rb\d v\\
&+\IInf g'(v)\lb\int_{0}^{v}\PPPp_\xi\left(\ln{\frac{v}{x}}\right)m(\d x)\rb\d v =:(g',F_{1}),
\end{split}
\]
where the interchange of integrals is permitted due to claims of Lemma \ref{auxLemma}.

Next, substituting $f'(x)=g(x)/x$ and $f''(x)=g'(x)/x-g(x)/x^{2}$, we get
\[
\begin{split}
\IInf \mathfrak{L}^{(\eta)}f(x) m(\d x)&=b_{\eta}\IInf \frac{g(x)}{x}m(\d x)+\frac{\sigma_{\eta}^{2}}{2}\IInf\bigg(\frac{g'(x)}{x}-\frac{g(x)}{x^{2}}\bigg)m(\d x)\\
&+\IInf \IInf \left(\frac{g'(x+w)}{x+w}-\frac{g(x+w)}{(x+w)^{2}}\right)\PPPp_{\eta}(w)\d wm(\d x)\\
&+\IInf \IInf \left(\frac{g'(x-w)}{x-w}-\frac{g(x-w)}{(x-w)^{2}}\right)\PPPn_{\eta}(w)\d wm(\d x).
\end{split}
\]
Again, using the identity $g(x)=\int_{0}^{x}g'(v)\d v$ and  the fact that $g$ is a function with compact support on $(0,\infty)$, we get after careful calculations and an appeal again to Lemma \ref{auxLemma} for interchange of integration
\[
\begin{split}
\IInf \mathfrak{L}^{(\eta)}f(x) m&(\d x)=b_{\eta}\IInf g'(v)\int_{v}^{\infty}\frac{m(\d x)}{x}\d v+\frac{\sigma_{\eta}^{2}}{2}\left(\IInf g'(x)\frac{m(\d x)}{x}-\IInf g'(v)\int_{v}^{\infty}\frac{m(\d x)}{x^{2}}\d v\right)\\
&+\IInf g'(v)\frac{1}{v}\int_{0}^{v}\PPPp_{\eta}(v-x)m(\d x)\d v-\IInf g'(w)\int_{w}^{\infty}\frac{1}{v^{2}}\int_{0}^{v}\PPPp_{\eta}(v-x)m(\d x)\d v\d w\\
&+\IInf g'(v)\frac{1}{v}\int_{v}^{\infty}\PPPn_{\eta}(x-v)m(\d x)\d v-\IInf g'(w)\int_{w}^{\infty}\frac{1}{v^{2}}\int_{v}^{\infty}\PPPn_{\eta}(x-v)m(\d x)\d v\d w\\
&:=(g',F_{2}).
\end{split}
\] 
We arrange the above expressions in the form $\int g'(x)\nu(\d x)$, where $\nu(\d x):=F_{1}(\d x)+F_{2}(\d x)$ is the same as in 
Lemma \ref{auxLemma}. From Lemma \ref{auxLemma}, we conclude that $\nu(\d x)$ defines a finite measure on every compact subset of $(0,\infty)$ and henceforth we consider it as a distribution in Schwartz's sense. Thus we get 
\[
0=\IInf \mathfrak{L}^{(U)}f(x) m(\d x)=(g',\nu)=(g,\nu')=(xf',\nu')=(f',x\nu')=\left(f,(x\nu')'\right),
\]
for each infinitely differentiable function $f$ with  compact support in $(0,\infty)$ and derivatives in the sense of Schwartz. Therefore using Schwartz theory of distributions for $\nu(dx)$, we get
that $x\nu'(\d x)=C\d x$ and therefore
\[
\nu(\d x)=\left(C\ln{x}+D\right)\d x.
\]

Next, we show that $C=D=0$. Note that from \eqref{limitNU} with $a=1$, we have $\lim_{z\to+\infty}z^{-1}\int_{1}^{z}\nu(\d v)=0$. Comparing this with 
\beqq
0=\lim_{z\to+\infty}z^{-1}\int_{1}^{z}(C\ln{x}+D)\d x=\lim_{z\to+\infty}(C\ln{z}-C+D)
\eeqq
we verify that $C=D=0$.
 Thus the proof of Theorem \ref{Theorem_main_integral_equation} is complete.
\qed
\vspace{0.5cm}

The next result is an almost immediate corollary of Theorem \ref{Prop1}, and in particular of formula \eqref{SS1}. See also  Corollary \ref{InfDiff} for a stronger result in a particular case when $\eta$ is a Brownian motion with drift.
\begin{corollary}\label{Corollary10}
Assume that condition \eqref{assumption_first_moments} is satisfied. If $\sigma^{2}_{\xi}+\sigma^{2}_{\eta}>0$ then $m(\d x)$ has a continuous density on $\r\setminus \{0\}$. 
\end{corollary}
\begin{proof}
The absolute continuity of $I(\xi,\eta)$ and boundedness of its derivative on compact subsets of $(0,\infty)$, when $\sigma^{2}_{\xi}+\sigma^{2}_{\eta}>0$ is immediate from \eqref{SS1}. Let $k(x)$ be the density of $m(\d x)$. To show the continuity of $k(x)$, we investigate all integral terms in \eqref{SS1}: all of them, except possibly the ones involving $\PPPp_{\xi}$ and $\PPPn_{\xi}$, are clearly continuous. Let us check continuity of these remaining two terms. Fix $v>0$ and $v/4>a>0$. Note that, for any real $h$ such that $|h|<v/4$, we have
\begin{align*}&\int_{0}^{v+h}\PPPp_{\xi}\lb\ln{\frac{v+h}{x}}\rb k(x)\d x=\int_{0}^{v+h-a}\PPPp_{\xi}\lb\ln{\frac{v+h}{x}}\rb k(x)\d x+\int_{v+h-a}^{v+h}\PPPp_{\xi}\lb\ln{\frac{v+h}{x}}\rb k(x)\d x.
\end{align*}
As $\PPPp_{\xi}$ is continuous and decreasing we verify the dominated convergence theorem applies, as $h\rightarrow 0$, by bounding $\PPPp_{\xi}$ in the first term and $k(x)$ in the second. This shows that all integral terms in \eqref{SS1} are continuous in $v$ and hence $k(v)$ is continuous. The computation for $\PPPn_{\xi}$ is the same whereas for $v<0$ we study $I(\xi,-\eta)$ with the same effect.
\end{proof}

\section{Exponential functionals with respect to Brownian motion with drift}\label{section_BM_exp_functionals}

In the next two sections, we study the special case when $\eta_t=\mu t+\sigma B_t$ is a Brownian motion with drift, so that the exponential functional is now defined as \begin{equation}\label{ExpFunBM}
I_{\mu,\sigma}:=\int_{0}^{\infty}e^{\xi_{t-}}(\mu \d t + \sigma \d B_{t}).
\end{equation}
We still work under Assumption \ref{assumption1}, note that the condition $\e[|\eta_1|]<\infty$ is clearly satisfied.  
From now on, we assume that $\sigma>0$, and in order to simplify notations we will write $\psi(s)=\psi_{\xi}(s)$. 
Note that formula \eqref{ExpFunBM} implies $I_{\mu,\sigma}\stackrel{d}{=} \sigma I_{\mu/\sigma,1}$, therefore it is sufficient to study the exponential functional with $\sigma=1$. 

The following three quantities will be very important in what follows
\beq\label{def_rho_hat_rho_and_theta}
\nonumber
\rho&:=&\sup\{z\ge 0:\, \e\left[ e^{z\xi_1} \right]<\infty \}, \\
\hat \rho&:=&\sup\{z\ge 0:\, \e\left[ e^{-z\xi_1} \right]<\infty \}, \\ \nonumber
\theta&:=&\sup\{z\ge 0:\, \e\left[ e^{z\xi_1} \right]\le 1 \}.
\eeq
In view of \eqref{SpecialFormExponents}, it is clear that
\beqq
\rho=\sup\left\{z\ge 0:\, \int_1^{\infty} e^{zx} \Pi_{\xi}(\d x)<\infty \right\}, \;\;\; 
\hat \rho=\sup\left\{z\ge 0:\, \int_{1}^{\infty} e^{zx} \Pi_{\xi}(-\d x)<\infty \right\}. 
\eeqq
Thus $\rho>0$ ($\hat \rho>0$) if and only if the measure $\Pi_{\xi}(\d x)$ has exponentially decaying positive (negative) tail. 
In this case the L\'evy-Khintchine formula \eqref{SpecialFormExponents} implies that the Laplace exponent $\psi(z)$ can be extended analytically in a strip $-\hat \rho < \re(z) < \rho$. It is clear from \eqref{def_rho_hat_rho_and_theta} that $0\le \theta\le \rho$. 
At the same time, due to Assumption \ref{assumption1} we have $\e[\xi_1]=\psi'(0)<0$ , which implies that 
$\theta>0$ if and only if $\rho>0$.

In the next Lemma we collect some simple analytical properties of the Laplace exponent $\psi(z)$. 
\begin{lemma}\label{BmLemma3}
Assume that $\xi$ satisfies condition \eqref{assumption_first_moments} and that $\rho>0$. Then $\psi(s)$ has no zeros in the strip $0<\re(s)<\theta$. Moreover if $\xi$ has a non-lattice distribution and $\psi(\theta)=0$, then $\theta$ is the unique zero of $\psi(s)$ in the strip $0<\re(s)\leq \theta$ and the unique real zero in the interval $(0,\rho)$.
\end{lemma}
\begin{proof}
Assume that $0<\re(s)<\theta$. Since 
\beqq
e^{\re(\psi(s))}=\left|\e \left[ e^{s\xi_{1}}\right]\right|\leq \e \left[ e^{\re(s)\xi_{1}}\right]=e^{\psi(\re(s))}
\eeqq
we conclude that $\re(\psi(s))\le\psi(\re(s))<0$, therefore $\psi(s)\ne 0$ in the strip $0<\re(s)<\theta$. 

Next, assume that $\psi(\theta+iy)=0$ for some $y\neq 0$ and $\xi$ has a non-lattice distribution. Then the characteristic function of the probability measure $e^{\theta v}\p (\xi_{1}\in dv)$ is equal to one at $y$, therefore it has to be a lattice distributed probability measure, see \cite[p 306, Theorem 5]{Shiryaev} which contradicts our assumption.

In order to prove that $\theta$ is the unique real zero of $\psi(s)$ on the interval $(0,\rho)$, we note that 
the first formula in \eqref{SpecialFormExponents} implies that
\[\psi''(s)=\sigma_{\xi}^{2}+\int_{\r}x^{2}e^{sx}\Pi_{\xi}(\d x)>0,\]
therefore $\psi(s)$ is convex on $(0,\rho)$ and it has at most one positive root at $\theta$.
\end{proof}

Next, let us introduce two other important objects
\beq\label{J}
J_{\alpha}:=\int_{0}^{\infty}e^{\alpha\xi_{t}}dt, \;\;\; {\textnormal{ and }} \;\;\;
V:=\frac{J_1^2}{J_2}.
\eeq
We will frequently use the following result, its proof follows immediately from Lemma 2.1 in \cite{MZ}: 
\begin{proposition}\label{Proposition1}
Assume that $\xi$ satisfies condition \eqref{assumption_first_moments}. For all $z\in \c$ in the strip $-1 \le \re(z)<\theta/\alpha$ we have $\e\left[J^{z}_{\alpha} \right]<\infty$.
\end{proposition}

Our main object of interest is the probability density function of $I_{\mu,\sigma}$, which we will denote by $k(x)$ (or by $k_{\mu,\sigma}(x)$ if we need to stress dependence on parameters). In the next Lemma, we collect some simple properties of $k(x)$.
\begin{lemma}\label{BmLemma1}
 Assume that $\xi$ satisfies condition \eqref{assumption_first_moments}. The law of $I_{\mu,\sigma}$ has a continuously differentiable density $k_{\mu,\sigma}(x)$ which is given by
\begin{equation}\label{BmLemma1-1}
k_{\mu,\sigma}(x)=\iint\limits_{\r_+^2}\frac{1}{\sigma\sqrt{2\pi z}}e^{-\frac{(x-\mu y)^{2}}{2z\sigma^{2}}}\p(J_{1}\in \d y;J_{2}\in \d z).
\end{equation}
Moreover, both functions $k_{\mu,\sigma}(x)$ and $k_{\mu,\sigma}'(x)$ are uniformly bounded on $\r$ and if $\mu\le 0$ then $k_{\mu,\sigma}(x)$ is decreasing on $\r_+$. 
\end{lemma}
\begin{proof}
Expression \eqref{BmLemma1-1} follows by  conditioning on $\xi$ and the fact that
\beqq
\int_{0}^{\infty}e^{f(t)}(\mu \d t + \sigma \d B_t) \stackrel{d}{=} N\left(\mu\int_{0}^{\infty}e^{f(s)}\d s;\sigma^{2}\int_{0}^{\infty}e^{2f(s)}\d s\right),
\eeqq
where $N(a,b)$ denotes a normal random variable with mean $a$ and variance $b$. 
The continuity of $k_{\mu,\sigma}(x)$ follows from the Dominated Convergence Theorem and the fact that $\e \left[J^{-\frac{1}{2}}_{2}\right]<\infty$, see Proposition \ref{Proposition1}.

Next, we observe that the function $|v|e^{-v^{2}}$ is bounded on $\mathbb{R}$ and therefore for some $C>0$ we have
\beqq
\bigg|\iint\limits_{\r_+^2}\frac{(x-\mu y)}{\sigma^{3}\sqrt{2\pi z^{3}}}e^{-\frac{(x-\mu y)^{2}}{2z\sigma^{2}}}\p(J_{1}\in \d y;J_{2}\in \d z)\bigg|\leq C \e [J^{-1}_{2}]<\infty,
\eeqq
where the last inequality follows from Proposition \ref{Proposition1}. This shows that we can differentiate the right-hand side 
of (\ref{BmLemma1-1}) and obtain
\begin{equation}\label{BmLemma1-2}
k'_{\mu,\sigma}(x)=-\iint\limits_{\r_+^2}\frac{(x-\mu y)}{\sigma^{3}\sqrt{2\pi z^{3}}}e^{-\frac{(x-\mu y)^{2}}{2z\sigma^{2}}}\p(J_{1}\in \d y;J_{2}\in \d z),
\end{equation}
and from the above discussion it follows that $|k_{\mu,\sigma}'(x)|\leq C \e [J^{-1}_{2}]<\infty$ for all $x\in \r$. Finally, for $\mu\leq 0$ and $x>0$ we check that $k_{\mu,\sigma}'(x)<0$ (see \eqref{BmLemma1-2}), therefore $k_{\mu,\sigma}(x)$ is decreasing.
\end{proof}

Our main tool for studying the properties of $k_{\mu,\sigma}(x)$ will be the Mellin transform of $I_{\mu,\sigma}$, which is defined for $\re(s)=1$ as
\beq\label{def_Mellin_transform_of_I}
\mathcal{M}_{\mu,\sigma}(s):=\e[(I_{\mu,\sigma})^{s-1}{\mathbf 1}_{\{I_{\mu,\sigma}>0\}}]=
\int_0^{\infty} x^{s-1} k_{\mu,\sigma}(x) \d x.
\eeq
Later we will extend this definition for a wider range of $s$, but a priori  it is not clear why this object should be finite 
for $\re(s)\ne 1$. 
Also, this choice of truncated random variable may seem awkward, since we only use the information about the density $k_{\mu,\sigma}(x)$ for $x\ge 0$. 
However, it is easy to see that the Mellin 
transform $\mathcal{M}_{\mu,\sigma}(s)$ uniquely determines $k_{\mu,\sigma}(x)$ for $x\ge 0$ while
$\mathcal{M}_{-\mu,\sigma}(s)$ uniquely determines $k_{\mu,\sigma}(x)$ for $x\le 0$. 
 This follows from the simple fact that $k_{\mu,\sigma}(-x)=k_{-\mu,\sigma}(x)$ (clearly $I_{\mu,\sigma}\stackrel{d}{=} -I_{-\mu,\sigma}$, see \eqref{ExpFunBM}). 
Moreover, later it will be clear that our definition of the Mellin transform is in fact quite natural, 
since $\mm_{\mu,\sigma}(s)$ satisfies the crucial functional equation \eqref{RecurBr}, 
which will lead to a wealth of interesting information about $k_{\mu,\sigma}(x)$.

As a first step in our study of the Mellin transform $\mm_{\mu,\sigma}(s)$ we obtain its analytic continuation  
 into a vertical strip in the complex plane. 
\begin{lemma}\label{Lemma_Mellin_transform_analytic_continuation}
 Assume that $\xi$ satisfies condition \eqref{assumption_first_moments}. The function ${\mathcal{M}}_{\mu,\sigma}(s)$ can be extended to an analytic function in the strip $-1<\re(s)<1+\theta$, except for a simple pole 
at $s=0$ with residue $k(0)$. Moreover, 
for all $s$ in the strip $-1<\re(s)<1+\theta$ we have
\beq\label{analytic_continuation_1}
\mathcal{M}_{\mu,\sigma}(s)=\frac{k(0)}{s}+\int_0^1 (k(x)-k(0)) x^{s-1} \d x+ \int_1^{\infty} k(x) x^{s-1} \d x,
\eeq
and for all $s$ in the strip $-1<\re(s)<0$ it is true that
\beq\label{analytic_continuation_2}
\mathcal{M}_{\mu,\sigma}(s)=-\frac{1}{s} \int_0^{\infty} x^s k'(x) \d x.
\eeq
\end{lemma}
\begin{proof}
First of all, since $k(x)$ is a probability density, it is integrable on $[0,\infty)$. Also, 
due to Lemma \ref{BmLemma1}, we know that $k(x)=k(0)+k'(0)x+o(x)$ as $x\to 0^+$, these two facts imply that $\mathcal{M}_{\mu,\sigma}(s)$ exists for all $s$ 
in the strip $0<\re(s)\le 1$.

Next, one can easily check that identity (\ref{analytic_continuation_1}) is valid for $s$ in the strip $0<\re(s)\le 1$. Since 
$k(x)-k(0)=k'(0)x+o(x)$, as $x\to 0^+$ we see that the first integral in the right-hand side of (\ref{analytic_continuation_1}) extends analytically
into the larger strip $-1<\re(s)<1$, while the second integral is analytic in the half-plane $\re(s)<1$. Thus 
(\ref{analytic_continuation_1}) provides an analytic continuation of $\mathcal{M}_{\mu,\sigma}(s)$ into the strip $-1<\re(s)<1$ and it is clear 
that $\mathcal{M}_{\mu,\sigma}(s)$ has a simple pole at $s=0$ with residue $k(0)$. 

Next, we note that for $-1<\re(s)<0$ we have
\beqq
\int_1^{\infty} k(x) x^{s-1} \d x=\int_1^{\infty} (k(x)-k(0)) x^{s-1} \d x-\frac{k(0)}{s}. 
\eeqq
Combining this expression with (\ref{analytic_continuation_1}) and applying integration by parts we obtain  (\ref{analytic_continuation_2}).

If $\theta=0$, then the proof is finished. However, if $\theta>0$  we still have to prove that $\mathcal{M}_{\mu,\sigma}(s)<\infty$ for $1<s<1+\theta$, 
and this requires a little bit more work. The proof will be based on certain special functions. The confluent hypergeometric function (see section 9.2 in \cite{Jeffrey2007} or chapter 6 in \cite{Erdelyi}) is defined as 
\beq\label{def_1F1}
{}_1F_1(a,b,z)=\sum\limits_{n\ge 0} \frac{(a)_n}{(b)_n} \frac{z^n}{n!},
\eeq
where $(a)_n=a(a+1)\dots(a+n-1)$ is the Pochhammer symbol. Using the ratio test it is easy to see that the series in 
\eqref{def_1F1} converges for all $z\in \c$, thus ${}_1F_1(a,b,z)$ is an entire function of $z$. We will also need the parabolic cylinder function, which is defined as
\beq\label{def_parabolic_cylinder}
D_p(z)=2^{\frac{p}2} e^{-\frac{z^2}4} \left[ \frac{\sqrt{\pi}}{\Gamma\left(\frac{1-p}2\right)} {}_1F_1\left(-\frac{p}2,\frac12; \frac{z^2}2\right)
+\frac{\sqrt{2 \pi} z}{\Gamma\left(-\frac{p}2\right)} {}_1F_1\left(\frac{1-p}2,\frac32; \frac{z^2}2\right)\right].
\eeq
Note that the parabolic cylinder function is analytic function of $p$ and $z$. 
See sections 9.24-9.25 in \cite{Jeffrey2007} for more information on the parabolic cylinder function. We will prove that 
$\mathcal{M}_{\mu,1}(s)$ exists for all $s$ in the strip $\re(s) \in (0,1+\theta)$ and everywhere in this strip we have
\beq\label{eqn_M_parabolic_cylinder}
\mathcal{M}_{\mu,1}(s)=\frac{\Gamma(s)}{\sqrt{2\pi}} \e\left[ J_2^{\frac{s-1}2} e^{-\frac{\mu^2}4 V} D_{-s}\left(-\mu \sqrt{V}\right) \right].
\eeq

Let us assume first that $\re(s)=1$. Then using (\ref{BmLemma1-1}) and (\ref{def_Mellin_transform_of_I}) we conclude that 
\beqq
\mathcal{M}_{\mu,1}(s)&=&\e \left[ \int_0^{\infty} x^{s-1} \frac{1}{\sqrt{2\pi J_2} }e^{-\frac{(x-\mu J_1)^2}{2J_2}} \d x \right]=\frac{1}{\sqrt{2\pi }} \e \left[ \frac{1}{\sqrt{J_2}} e^{-\frac{\mu^2}2 V} \int_0^{\infty} x^{s-1} e^{-\frac{1}{2J_2} x^2 + \frac{\mu J_1}{J_2} x} \d x \right].
\eeqq
Performing the change of variables $x=u \sqrt{J_2}$ and using the following integral identity (formula 9.241.2 in \cite{Jeffrey2007})
\beqq
\int_{0}^{\infty} u^{s-1} e^{-\frac{u^2}2-uz} \d u = \Gamma(s) e^{\frac{z^2}4} D_{-s}(z), \;\;\; \re(s)>0,
\eeqq
we obtain equation (\ref{eqn_M_parabolic_cylinder}).

Thus, we have established that (\ref{eqn_M_parabolic_cylinder}) is true for all $s$ on the vertical line $\re(s)=1$. Now, we will perform analytic continuation into a larger domain. Formulas 9.246 in  \cite{Jeffrey2007}, give us the following asymptotic expansions: for $z\in \r$ 
\beq\label{Dsz_asymptotics}
D_{-s}(z)=
\begin{cases}
 O\left(z^{-s} e^{-\frac{z^2}4} \right), \qquad\qquad \qquad\;\;\;\;\; z\to +\infty,\\ 
 O\left(z^{-s} e^{-\frac{z^2}4} \right)+O\left(z^{s-1} e^{\frac{z^2}4} \right), \;\;\; z\to -\infty.
\end{cases}
\eeq
Assume that $\mu<0$ and $s\in (0,1+\theta)$ or $\mu>0$ and $s\in (0,1)$. Then, from  (\ref{Dsz_asymptotics}) and the fact that $D_s(z)$ is a continuous
function of $z$ we find that there
exists a constant $C_1>0$ such that $|e^{-\frac{\mu^2}4 z} D_{-s}\left(-\mu \sqrt{z}\right)|<C_1$ for all $z>0$. Therefore from
 (\ref{eqn_M_parabolic_cylinder}), we conclude that
\beqq
\vert \mathcal{M}_{\mu,1}(s)\vert < C_1  \frac{\Gamma(s)}{\sqrt{2\pi}} \e\left[ J_2^{\frac{s-1}2}\right],
\eeqq
and the right-hand side is finite if $s\in (0,1+\theta)$, see Proposition \ref{Proposition1}.

Next, when $\mu>0$ and $s \in (1, 1+\theta)$, we again use  (\ref{Dsz_asymptotics}) and the fact that $D_s(z)$ is continuous in $z$ 
to conclude that there exists $C_2>0$ such that $|e^{-\frac{\mu^2}4 z} D_{-s}\left(-\mu \sqrt{z}\right)|<C_2 z^{\frac{s-1}2}$ for all $z>0$.
Therefore from
 (\ref{eqn_M_parabolic_cylinder}), we conclude that
\beqq
\vert \mathcal{M}_{\mu,1}(s)\vert < C_2  \frac{\Gamma(s)}{\sqrt{2\pi}} \e\left[ J_2^{\frac{s-1}2} V^{\frac{s-1}2}\right]
=C_2  \frac{\Gamma(s)}{\sqrt{2\pi}} \e\left[ J_1^{s-1}\right],
\eeqq
and the right-hand side is finite if $s\in (1,1+\theta)$, see Proposition \ref{Proposition1}. 
\end{proof}

The next theorem is our first main result in this section. 
\begin{theorem}\label{BmTheorem1} Assume that $\xi$ satisfies condition \eqref{assumption_first_moments} and that  $\theta>0$. Then for all $s$ such that $0<\re(s)<\theta$, we have
\beq\label{RecurBr}
\frac{\psi(s)}{s}{\mathcal{M}}_{\mu,\sigma}(s+1)+\mu \mm_{\mu,\sigma}(s)+ \frac{\sigma^{2}}{2}(s-1) \mm_{\mu,\sigma}(s-1)=0.
\eeq
\end{theorem}
\begin{proof}
Setting $\Pi_{\eta}\equiv 0$ in \eqref{SS1} we find that $k(x)$ satisfies the following integral equation
\beq\label{eqn_kx_n1}
\frac{\sigma_{\xi}^{2}}{2} F_1(k;v)+b_{\xi} F_2(k;v)+
F_3(k;v)+F_4(k;v)+\mu F_5(k;v)+\frac{\sigma^2}2 F_6(k;v)=0,
\eeq
where we have defined $F_1(k;v):=k(v)$ and
\beq\label{def_Fi}
\nonumber
&&F_2(k;v):=\frac{1}{v}\int_{v}^{\infty}k(x)\d x, \;\;\;
F_3(k;v):=\frac{1}{v}\int_{v}^{\infty}\PPPn_{\xi}\left(\ln\left(\frac{x}{v}\right)\right)k(x)\d x, \\
&&
F_4(k;v):=\frac{1}{v}\int_{0}^{v}\PPPp_{\xi}\left(\ln\left(\frac{v}{x}\right)\right)k(x)\d x, \;\;\;
F_5(k;v):=\frac{1}{v}\int_{v}^{\infty}\frac{k(x)}{x}\d x, \\ \nonumber
&&
F_6(k;v):=\frac{1}{v} \left[\frac{k(v)}{v}-\int_{v}^{\infty}\frac{k(x)}{x^{2}}\d x \right].
\eeq
Our plan is to compute the Mellin transform of each term in (\ref{eqn_kx_n1}). Assume first that $1<\re(s)<1+\min\{1,\theta\}$. According to Lemma 
\ref{Lemma_Mellin_transform_analytic_continuation},
the Mellin transform of the first term exists in this strip and is equal to $\mathcal{M}_{\mu,\sigma}(s)$.

Let us compute the Mellin transform of the second term. We use integration by parts and obtain for all $y>0$
\beq\label{proof_F2kv}
\int_0^{y} v^{s-1} F_2(k;v) \d v=\frac{1}{s-1} y^{s-1} \int_y^{\infty} k(x) \d x+\frac{1}{s-1} \int_0^y v^{s-1} k(v) \d v.
\eeq
As $y\to +\infty$ the first term in the right-hand side of the above equation goes to zero (this follows from the fact that the $k(x) x^{s-1}$
is absolutely integrable on $(0,\infty)$), thus we conclude that the Mellin transform of $F_2(k;v)$ is equal to $\mathcal{M}_{\mu,\sigma}(s)/(s-1)$. In  exactly the same way one finds that the Mellin transform of $F_5(k;v)$ is equal to $\mathcal{M}_{\mu,\sigma}(s-1)/(s-1)$.

Let us consider the third term $F_3(k;v)$. Performing the change of variables $x \mapsto yv$ we find that
\beqq
F_3(k;v)=\int_1^{\infty}\PPPn_{\xi}\left(\ln(y)\right)k(yv)\d y.
\eeqq
Therefore the Mellin transform of $F_3(k;v)$ is given by
\beqq
&&\int_0^{\infty} v^{s-1} F_3(k;v) \d v=\int_1^{\infty}\PPPn_{\xi}\left(\ln(y)\right) \int_0^{\infty} v^{s-1} k(yv) \d v \d y\\
&&=\left[\int_1^{\infty}\PPPn_{\xi}\left(\ln(y)\right) y^{-s} \d y \right] \times \mathcal{M}_{\mu,\sigma}(s)=
\left[\int_0^{\infty}\PPPn_{\xi}(u) e^{-(s-1)u} \d u \right] \times \mathcal{M}_{\mu,\sigma}(s),
\eeqq
where we have used Fubini's theorem in the first step, performed the change of variables $v \mapsto z/y$ in the second step and applied the 
change of variables $y \mapsto \exp(u)$ in the last step.

In exactly the same way we find that the Mellin transform of $F_4(k;v)$ is equal to
\beqq
\int_0^{\infty} v^{s-1} F_4(k;v) \d v=\left[\int_0^{\infty}\PPPp_{\xi}(u) e^{(s-1)u} \d u \right] \times \mathcal{M}_{\mu,\sigma}(s).
\eeqq

Finally, let us consider the sixth term $F_6(k;v)$. Using integration by parts and the fact that $k(x)$ is bounded we find that
\beq\label{formula_F6}
F_6(k;v)=-\frac{1}{v} \int_v^{\infty} \frac{k'(x)}{x} \d x.
\eeq
Since $k'(x)$ is uniformly bounded on $[0,\infty)$ we conclude that $F_6(k;v)=O(\ln(v)/v)$, as $v\to 0^+$, and from (\ref{def_Fi}) we see that
$F_6(k;v)=O(1/v^2)$, as $v\to +\infty$. This shows that the Mellin transform of $F_6(k;v)$ exists for $1<\re(s)<1+\min\{1,\theta\}$. Using
(\ref{formula_F6}) and integration by parts we find that for $0<v_0<v_1<\infty$
\beq\label{formula_F6_n2}
\int_{v_0}^{v_1} v^{s-1} F_6(k;v) \d v=\frac{v_1^{s}}{s-1} F_6(k;v_1) -
\frac{v_0^{s}}{s-1}  F_6(k;v_0) -\frac{1}{s-1}
\int_{v_0}^{v_1} v^{s-2} k'(v) \d v.
\eeq
From the above discussion we find that the first (second) term in right-hand side of (\ref{formula_F6_n2}) converges to zero as $v_1 \to +\infty$
($v_0 \to 0^+$), therefore from (\ref{analytic_continuation_2}) and (\ref{formula_F6_n2}) we conclude that for $1<\re(s)<1+\min\{1,\theta\}$ the Mellin transform of $F_6(k;v)$ is given by
\beqq
\int_{0}^{\infty} v^{s-1} F_6(k;v) \d v=-\frac{1}{s-1} \int_0^{\infty} v^{s-2} k'(v) d v=\frac{s-2}{s-1}\mathcal{M}_{\mu,\sigma}(s-2).
\eeqq

Collecting all the terms in (\ref{eqn_kx_n1}) we see that for all $s$ in the strip $1<\re(s)<1+\min\{1,\theta\}$ we have
\beqq
\frac{\sigma_{\xi}^2}{2} \mathcal{M}_{\mu,\sigma}(s)&+& \frac{b_{\xi}}{s-1}\mathcal{M}_{\mu,\sigma}(s)+
\left[\int_0^{\infty}\PPPn_{\xi}(u) e^{-(s-1)u} \d u+\int_0^{\infty}\PPPp_{\xi}(u) e^{(s-1)u} \d u \right] \times \mathcal{M}_{\mu,\sigma}(s)\\
&+&\frac{\mu}{s-1} \mathcal{M}_{\mu,\sigma}(s-1)+\frac{\sigma^2}{2}\frac{s-2}{s-1}\mathcal{M}_{\mu,\sigma}(s-2)=0.
\eeqq
Formula (\ref{RecurBr}) follows from the above equation by changing variables $s\mapsto s-1$ and applying formula
(\ref{SpecialFormExponents}). This ends the proof in the case $\theta \in (0,1)$. If $\theta>1$ then 
 \eqref{RecurBr} can be extended from the strip $0<\re(s)<\min\{1,\theta\}$ to $0<\re(s)<\theta$ by analytic continuation. 
\end{proof}

\begin{remark}\label{Th2Rem1}
Note that the functional equation \eqref{RecurBr} is a more general version of the well-known functional equation 
when $\sigma=0$, see formula (2.3) in Maulik and Zwart, \cite{MZ}. Nonetheless, the derivation of \eqref{RecurBr} requires 
the integral equation \eqref{SS1} whereas the classical functional equation (2.3) in \cite{MZ} can be obtained by rather simple arguments.
\end{remark}

Theorem \ref{BmTheorem1} will prove crucial for applications. It allows to derive the analytical properties of the Mellin transform 
$\mm_{\mu,\sigma}(s)$ (such as its behaviour at the singularities and their precise location in the complex plane) from the  properties of the Laplace exponent $\psi(s)$ itself. The next result serves to illustrate these ideas.

\begin{corollary}\label{Corollary11} Assume that $\xi$ satisfies condition \eqref{assumption_first_moments} and that $\theta>0$. 
\\
\begin{itemize}
 \item[(i)]  The function $\mathcal{M}_{\mu,\sigma}(s)$ can be analytically continued into the strip
$\re(s)\in (-1-\hat \rho,1+\rho)$. Its only singularities in the strip $-1-\hat \rho < \re(s) < 1+\theta$ are the simple poles at the points $\{-n \; : \; 0\le n < 1+\hat \rho\}$. 
 \item[(ii)] If $\xi$ has a non-lattice distribution and $\theta<\rho$ then $\mathcal{M}_{\mu,\sigma}(s)$ has a simple pole at $s=1+\theta$ with 
residue
\beq\label{def_C}
R(\theta):=-\frac{\theta}{\psi'(\theta)}\lb \mu\mathcal{M}_{\mu,\sigma}(\theta)+\frac{\sigma^{2}}2(\theta-1)\mathcal{M}_{\mu,\sigma}(\theta-1) \rb.
\eeq
 The only other singularities of $\mathcal{M}_{\mu,\sigma}(s)$ in the strip 
$0<\re(s) < 1+\rho$ are poles of the form $\zeta+n$, where $n\in {\mathbb{N}}$ and $\zeta$ is 
a root of $\psi(s)$ in the strip $\theta<\re(s)<\rho$.
\item[(iii)] Consider the ``boundary'' case $\theta=\rho$. Assume that $\xi$ has a non-lattice distribution. If $\psi(\theta)<0$ 
 then the function $\mathcal{M}_{\mu,\sigma}(s)$ is continuous in the strip $0<\re(s)\le 1+\theta$. On the other hand, if $\psi(\theta)=0$  and 
$\e\left[ \xi_1^2 \exp(\theta \xi_1)\right]<\infty$, then the function $\mathcal{M}_{\mu,\sigma}(s)-R(\theta)/(s-1-\theta)$ is continuous in the strip $0<\re(s)\le 1+\theta$. 
\end{itemize}
\end{corollary}

\begin{proof}
The proof of parts (i) and (ii) follows easily from Theorem \ref{BmTheorem1} and Lemmas \ref{BmLemma3} and \ref{Lemma_Mellin_transform_analytic_continuation}.

Let us prove (iii). If $\psi(\theta)<0$ we use the same argument as in the proof of Lemma \ref{BmLemma3} and conclude that 
$\re(\psi(s))\le\psi(\re(s))=\psi(\theta)<0$ on the line $\re(s)=\theta$; this fact and Lemma \ref{BmLemma3} imply that $\psi(s)\ne 0$ in the strip $0<\re(s)\le \theta$. Since $\psi(\theta)<0$ we can use \eqref{SpecialFormExponents} and the dominated convergence theorem 
to show that $\psi(s)$ is continuous in the strip $0<\re(s)\le \theta$. These two facts and the functional equation \eqref{RecurBr}
show that $\mm_{\mu,\sigma}(s)$ is continuous in the strip $0<\re(s)\le \theta$. 

Finally, let us consider the case when $\theta=\rho$ and $\psi(\theta)=0$. Condition $\e[ \xi_1^2 \exp(\theta \xi_1)]<\infty$ and the dominated 
convergence theorem show that the functions $\psi(s)$, $\psi'(s)$ and $\psi''(s)$, which are analytic in the strip $0<\re(s)<\theta$, can be continuously extended to the right boundary of this strip $\re(s)=\theta$.
 Again, using \eqref{SpecialFormExponents} and the dominated convergence theorem one can check that as $s\to \theta$ in the strip $0<\re(s)\le \theta$, it is true that
\beq\label{substract_pole}
H(s):=\frac{1}{\psi(s)}-\frac{1}{\psi'(\theta)(s-\theta)}\to-\frac12 \frac{\psi''(\theta)}{\psi'(\theta)^2}<\infty. 
\eeq
Note that $\psi'(\theta)>0$ due to the convexity of $\psi(s)$ on the interval $0<s<\rho$. 
Lemma \ref{BmLemma3} and the fact that $\xi$ has non-lattice distribution guarantee that the only zero of $\psi(s)$ in the strip $0<\re(s)\le \theta$ 
is at $s=\theta$. From here and from \eqref{substract_pole}, we see that the function $H(s)$ defined in \eqref{substract_pole}
is continuous in the strip $0<\re(s)\le \theta$. 

Let us define 
\beqq
F(s)=-s\left( \mu \mm_{\mu,\sigma}(s) + \frac{\sigma^2}2 (s-1) \mm_{\mu,\sigma}(s-1) \right).
\eeqq
It is clear from Lemma \ref{Lemma_Mellin_transform_analytic_continuation} that $F(s)$ is analytic in some neighbourhood of the line $\re(s)=\theta$, thus
the function  $G(s)=(F(s)-F(\theta))/(s-\theta)$ is also analytic in the neighbourhood of the line $\re(s)=\theta$. 
Next, we use the functional equation \eqref{RecurBr} in the form $\mm_{\mu,\sigma}(s+1)=F(s)/\psi(s)$ and after rearranging the terms, we find 
\beqq
\mm_{\mu,\sigma}(s+1)-\frac{F(\theta)}{\psi'(\theta)} \frac{1}{s-\theta}= F(s) 
H(s) + \frac{G(s)}{\psi'(\theta)}.
\eeqq
From the above discussion it is clear that the function in the right-hand side is continuous in the strip $0< \re(s) \le \theta$, which ends the proof of part (iii).
\end{proof}

In view of \eqref{BmLemma1-1} it is clear that $k_{\mu,\sigma}(x)$ depends on the joint distribution of $J_1$ and $J_2$. As we will see later 
in Proposition \ref{AppProp2}, the Mellin transform $\mm_{\mu,\sigma}(x)$ can be expressed in terms of the joint moments $\e [J_1^u J_2^v]$. The next Lemma
presents several crucial results on the existence of joint moments of this form. Recall that  $V=J_1^2/J_2$.

\begin{lemma}\label{ProbaLemma} Assume that $\xi$ satisfies condition \eqref{assumption_first_moments}. 
${}$
\\
\begin{itemize}
 \item[(i)] There exists $\epsilon>0$ such that
\begin{equation}\label{FiniteMean}
\e \left[e^{\epsilon V} \right]<\infty.
\end{equation}
 \item[(ii)]
For any $(u,s)\in \r^2$ in the domain 
\beqq
\mathcal{D}=\Big\{-1<s<1+\theta,\, u\leq 0\Big\}\cup\Big\{s>0; u>0;\,u\leq 1-s\Big\}
\eeqq
we have
\begin{equation}\label{FinitenessMoments}
\e\left[J^{-u}_{1}J^{\frac{1}{2}(u+s-1)}_{2}\right]<\infty.
\end{equation}
The function $(u,s)\in \c^2 \mapsto \e\left[J^{-u}_{1}J^{\frac{1}{2}(u+s-1)}_{2}\right]$ is analytic as long as $(\re(s),\re(u))\in {\mathcal{D}}$ and it
is uniformly bounded if $(\re(s),\re(u))$ belongs to a compact subset of ${\mathcal{D}}$.
\end{itemize}
\end{lemma}

\begin{proof}
Let us prove (i). Denote by $J_{1}(t)=\int_{0}^{t}e^{\xi_{s}}\d s$ and $J_{2}(t)=\int_{0}^{t}e^{2\xi_{s}}\d s$. 
It is clear that $J_1(0)=J_2(0)=0$ and that both $J_1(t)$ and $J_2(t)$ are continuous in $t$. Since 
\beq\label{der_Ji}
\frac{\d }{\d t} J_1^2(t) \bigg|_{t=0}=0 \qquad {\textnormal{ and }}\qquad \frac{\d }{\d t} J_2(t) \bigg|_{t=0}=1, 
\eeq 
we conclude that for every $x>0$ with probability one, there exists $\epsilon>0$ such that $J_1^2(t)<x J_2(t)$ for $0<t<\epsilon$. 
This fact and the continuity of $J_1(t)$ and $J_2(t)$ imply that
\beq\label{def_gx}
g(x):=\p\left(T_x<\infty\right) \ge \p(J^{2}_{1}>xJ_{2}),
\eeq
where  $T_x=\inf\{t>0\,:\,J^{2}_{1}(t)/J_2(t)=x\}$ and as usual we assume that $\inf\{\emptyset\}=+\infty$.
 We aim to show that for all $x>0$, we have $g(2x)\leq g^{2}(x)$.

 From the (\ref{der_Ji}), we know that
 $J_1^2(t)/J_2(t)\to 0$ as $t\to 0^+$. This fact and the continuity of $J_1(t)$ and 
$J_2(t)$ imply that $\p(T_x=0)=0$ for all $x>0$ and $T_x<T_y$ a.s., for $y>x$. Using the inequality $2a^{2}+2b^{2}\geq (a+b)^{2}$, we get 
\beqq
2\left(J_{1}(T_x+t)-J_{1}(T_x)\right)^{2}+2J^{2}_{1}(T_x)\geq J^{2}_{1}(T_x+t),
\eeqq
 and we estimate
\begin{align*}
&g(2x)=\p(T_{2x}<\infty)=\p\left(T_x<\infty;\exists\,t>0\,: J_{1}(T_x+t)^{2}=2x J_2(T_x+t) \right)\leq \\
&\p\Big(T_x<\infty;\exists\,t>0\,: 2\left(J_{1}(T_x+t)-J_{1}(T_x)\right)^{2}+2J^{2}_{1}(T_x)=2x J_2(T_x+t)\Big).
\end{align*}
Since $J^{2}_{1}(T_x)=xJ_{2}(T_x)$, we obtain from the above inequality
\begin{align*}
&g(2x)\leq \p\left(T_x<\infty;\exists\,t> 0\,:(J_{1}(T_x+t)-J_{1}(T_x))^{2}=x(J_{2}(T_x+t)-J_{2}(T_x))\right)=\\
& \p\left(T_x<\infty;\exists\,t> 0\,:e^{2\xi_{T_x}}\tilde{J}^{2}_{1}(t)=xe^{2\xi_{T_x}}\tilde{J}_{2}(t)\right)=g^{2}(x),
\end{align*}
where $\tilde{J}_{i}$ are the exponential functionals based on $\tilde\xi_t=\xi_{T_x+t}-\xi_{T_x}$ and we have used the fact
that the process $\{\tilde\xi_t\}_{t\ge 0}$ is independent of ${\mathcal F}_{T_x}$. Thus, we have obtained the key inequality 
$g(2x)\le g^2(x)$.

Next, let us prove that there exists $x^*>0$ such that $g(x^*)<1$. Assume that the converse is true, that is $g(x)=1$ for all $x>0$. 
In particular, $g(n)=1$ for all $n\geq 1$. Let $A_n=\{\exists \,t>0\, : J_1(t)^2=nJ_2(t)\}$. Since $\p(A_n)=1$ for all $n\ge 1$, we conclude
that $\p(\cap_{n\ge 1} A_n)=1$. This implies that with probability one there exists a strictly increasing random sequence of positive numbers $\{t_{n}\}_{n\geq 1}$ such that $J^{2}_{1}(t_{n})=nJ_{2}(t_{n})$. Since $J_2(t_n)$ is an increasing sequence, we conclude that as $n\to +\infty$ we have
$\p(J_1(t_n)^2 \to +\infty)=1$, and due to the fact that $J_1(t_n)\le J_1$,  we arrive at a contradiction $\p(J_1=\infty)=1$. 

Thus, we have proved that there exists $x^*>0$ such that $g(x^*)<1$.  For $x>x^*$, let us define $N>0$ to be the unique integer number such that $2^N \le x/x^*< 2^{N+1}$. Applying the
inequality $g(x)<g(x/2)^2$ exactly $N$ times, we obtain $g(x) \le g(x/2^N)^{2^N}$. Using the fact that $g(x)$ is a decreasing function and that
 $x^* \le x/2^N$, we conclude that $g(x)\le g(x^*)^{2^N}$, and since $x/(2x^*)<2^N$, we see that for all $x>x^*$ 
\beqq
g(x) < e^{-cx},
\eeqq 
where $c=- \ln (g(x^*))/(2x^*)>0$. This fact and (\ref{def_gx}) imply that  $\p(V > x)<\exp(-c x)$ for all $x> x^*$, 
thus (\ref{FiniteMean}) is true for any $\epsilon \in (0,c)$. This ends the proof of part (i).

Let us prove (ii).  Assume first that $u\leq 0$. Then using Holder inequality we get
\[\e\left[J^{-u}_{1}J^{\frac{1}{2}(u+s-1)}_{2}\right]=\e\left[V^{-\frac{u}{2}}J^{\frac{1}{2}(s-1)}_{2}\right]\leq \lb\e\left[V^{-\frac{u}{2}p}\right]\rb^{\frac{1}{p}}\lb\e\left[J^{\frac{q}{2}(s-1)}_{2}\right]\rb^{\frac{1}{q}}.\]
From part (i), we know that $V$ has finite positive moments of all orders. Then it suffices to choose $q=2(1-s)^{-1}$ for $-1<s<0$,  $q=2$ for $0\leq s\leq 1$ and $q=\frac{1}{2}+\frac{\theta}{2(s-1)}$ for $1<s<1+\theta$ to conclude that \eqref{FinitenessMoments} holds.

Assume next $0<s< 1$, $u>0$ and $u\leq 1-s$. Then with $p=u^{-1}$ and $q=(1-u)^{-1}$ we have 
\[
\e\left[J^{-u}_{1}J^{\frac{1}{2}(u+s-1)}_{2}\right]\leq \lb\e\left[ J^{-1}_{1}\right]\rb^{u}\lb\e\left[J^{\frac{1}{2}(-1+\frac{s}{1-u})}_{2}\right]\rb^{1-u}<\infty,\]
due to Proposition \ref{Proposition1} and the fact that $\lb-1+\frac{s}{1-u}\rb\in(-1,0]$. 
\end{proof}

Now we are ready to present several integral expressions for the Mellin transform $\mm_{\mu,\sigma}(s)$. These expressions are interesting in their own right, but they will also lead to an important result about the exponential decay of $\mm_{\mu,\sigma}(s)$ as $\im(s)\to \infty$ (Theorem 
\ref{thm_asymptotic_decay} below). Note that due to the identity $I_{\mu,\sigma}\stackrel{d}{=}\sigma I_{\mu/\sigma,1}$ we have $\mm_{\mu,\sigma}(s)=\sigma^{s-1} \mm_{\mu/\sigma,1}(s)$, therefore it is enough to state the results for $\sigma=1$.

\begin{proposition}\label{AppProp2} Assume that $\xi$ satisfies condition \eqref{assumption_first_moments}. 
${}$
\\
\begin{itemize}
\item[(i)]  For $-1<\re(s)<1+\theta$
\begin{equation}\label{AppLemma2-2}
\mm_{0,1}(s)=\frac{2^{-\frac12(s+1)}\Gamma(s)}{\Gamma\lb\frac12(s+1)\rb }\e \left[J_2^{\frac12 (s-1)} \right].
\end{equation}
\item [(ii)] For $\mu<0$ and $-1<\re(s)<1+\theta$
\beq\label{AppLemma2-1}
\mm_{\mu,1}(s)=\mm_{0,1}(s)+
\frac{2^{-\frac12(s+1)}}{2\pi \i} \int_{-\frac1{2}+\i \r} \frac{ \Gamma(s)\Gamma(u)}{\Gamma\lb \frac12(u+s+1)\rb} 
 \e\left[J_1^{-u}J_2^{\frac12(u+s-1)}\right] 
 (2\mu^2)^{-\frac{u}2} \d u.
\eeq
\item[(iii)] For $\mu>0$ and $-1<\re(s)<1+\theta$  
\begin{equation}\label{AppLemma2-3}
\mm_{\mu,1}(s)= \frac{2^{-\frac12(s-1)}\Gamma(s)}{\Gamma\lb\frac12(s+1)\rb }
\e \left[ J_2^{\frac12 (s-1)} {}_1F_1\left(\frac{1-s}2,\frac12,-\frac{\mu^2}2 V \right) \right]-\mm_{-\mu,1}(s).
\end{equation}
\end{itemize}
\end{proposition}

The proof of this Proposition is quite technical, therefore we have divided it into several steps. First of all, in 
Lemmas \ref{Gamma_estimates}, \ref{AppLemma1}, \ref{LastLemma}, we 
establish several technical results which will be needed in the proof of Proposition \ref{AppProp2}, and also useful  later.

\begin{lemma}\label{Gamma_estimates}

${}$ \\
\begin{itemize}
\item[(i)] For every $\epsilon>0$ and $a<b$ there exists $C=C(\epsilon, a,b)>0$ such that 
\beqq
\big | \Gamma(x+\i y) \big | < C e^{-\left( \frac{\pi}2-\epsilon\right) |y|}
\eeqq
for all $a<x<b$ and $|y|>1$. 
\item[(ii)] For every $\epsilon>0$ and $a<b$ there exists $C=C(\epsilon, a,b)>0$ such that 
\beqq
\big | \Gamma(x+\i y) \big |> C e^{-\left( \frac{\pi}2+\epsilon\right) |y|}
\eeqq
for all $a<x<b$ and $y\in \r$. 
\end{itemize}
\end{lemma}

\begin{proof}
We start with the following asymptotic expression
\beq\label{GammaAsympt}
\Big|\Gamma(x+iy)\Big|=\sqrt{2 \pi}|y|^{x-\frac{1}{2}}e^{-\frac{\pi}{2}|y|}\left(1+O\left(\frac{1}{|y|}\right)\right), \;\;\; y\to \infty
\eeq
which holds uniformly in $x$ on compact subsets of $\r$, see formula 8.328.1 in \cite{Jeffrey2007}.  Part (i) follows easily from 
(\ref{GammaAsympt}) and for part (ii), we use the additional fact that $\Gamma(s)$ has no zeros in the entire complex plane. 
\end{proof}

\begin{lemma}\label{AppLemma1}
For $\mu<0$, $\re(w)<\frac{1}{2}$ and $0<\re(s)<1-2\re(w)$
\begin{equation}\label{AppLemma1-1}
\iint_{\mathbb{R}^{2}_{+}}\frac{1}{\sqrt{2\pi z}}e^{-\frac{(x-\mu y)^{2}}{2z}}x^{s-1}z^{w-1}\d x\d z=2^{w-1}(-\mu y)^{2w+s-1}\frac{\Gamma(s)\Gamma(1-2w-s)}{\Gamma(1-w)}.
\end{equation}
\end{lemma}
\begin{proof}
We change the variable of integration $z\mapsto \frac{1}{u}$ and find that for $a>0$ and $\re(w)<1/2$
\[\int_{0}^{\infty}e^{-\frac{a}{z}}z^{w-1-\frac{1}{2}}dz=\int_{0}^{\infty}e^{-au}u^{-\frac{1}{2}-w}du=a^{w-\frac{1}{2}}\Gamma\Big(\frac{1}{2}-w\Big).\]
Then for $2\re(w)+\re(s)-1<0$ we can apply the Fubini's theorem and obtain
\beqq
&&\iint_{\mathbb{R}^{2}_{+}}\frac{1}{\sqrt{2\pi z}}e^{-\frac{(x-\mu y)^{2}}{2z}}x^{s-1}z^{w-1}\d x\d z=\frac{1}{\sqrt{2\pi}}2^{\frac{1}{2}-w}\Gamma\Big(\frac{1}{2}-w\Big)\int_{0}^{\infty}(x-\mu y)^{2w-1}x^{s-1}\d x\\
&&\qquad\qquad\qquad=(-\mu y)^{2w+s-1}\frac{1}{\sqrt{2\pi}}2^{\frac{1}{2}-w}\Gamma\Big(\frac{1}{2}-w\Big)\int_{0}^{\infty}(x+1)^{2w-1}x^{s-1}\d x\\
&&\qquad\qquad\qquad=(-\mu y)^{2w+s-1}\frac{1}{\sqrt{2\pi}}2^{\frac{1}{2}-w}\Gamma\Big(\frac{1}{2}-w\Big)
\frac{\Gamma(s)\Gamma(1-2w-s)}{\Gamma(1-2w)},
\eeqq
where in the last step we have used the beta-integral identity (see equation 3.194.3 \cite{Jeffrey2007}).
Formula \eqref{AppLemma1-1} can be derived from the above equation by application of the Legendre duplication formula for the gamma function (see formula 8.335.1 in \cite{Jeffrey2007}).
\end{proof}

\begin{lemma}\label{LastLemma}
Assume that $a_{0}<a_{1}$ and $b\in \c$ are such that $\re(b)\in (0,1)\cup(1,2)$. 
Recall that ${}_1F_1(a,b,z)$ denotes the confluent hypergeometric function defined by \eqref{def_1F1}. For each $\epsilon>0$ there exist 
a constant $C=C(a_{0},a_{1},b,\epsilon)>0$ and a constant $D=D(a_{0},a_{1},\epsilon)\in(0,\frac{\pi}{2})$ such that for all $a\in \mathbb{C}$ with $a_{0}<\re(a)<a_{1}$ and all $z>0$
\begin{equation}\label{ConluentAst}
\big|{}_1F_{1}(a,b,-z)\big|\leq C e^{\epsilon z+D|Im(a)|}.
\end{equation}
\end{lemma}
\begin{proof}
We start with the following integral representation
\begin{equation}\label{LastLemma1}
{}_1F_{1}(a,b,-z)=\frac{1}{2\pi i}\Gamma(b)z^{1-b}e^{-z}\int_{\gamma+i\mathbb{R}}e^{wz}w^{-b}(1-w^{-1})^{a-b}\d w,
\end{equation}
which holds for $z>0$, $\re(b)>0$ and $\gamma>1$. This representation follows from formula (7) on page 273 in \cite{Erdelyi} and 
the identity ${}_1F_1(a,b,-z)=\exp(-z) {}_1F_1(b-a,b,z)$ (see formula (7) on page 253 in \cite{Erdelyi}). 

Next fix $\epsilon>0$ and assume that $\re(b)\in (1,2)$ and $z\geq 1$. We also denote $\gamma=1+\epsilon$. Then
changing variables $w \mapsto \gamma + \i t$ we obtain from \eqref{LastLemma1}
\beq\label{LastLemma3}
\nonumber
\big|{}_1F_{1}(a,b,z)\big|&=&\Bigg|\frac{1}{2\pi i}\Gamma(b)z^{1-b}e^{\epsilon z}\int_{-\infty}^{\infty}e^{itz}(\gamma+it)^{-b}\Big(1-(\gamma+it)^{-1}\Big)^{a-b}dt\Bigg|\\
&\leq&
 C_1(b)e^{\epsilon z}\int_{-\infty}^{\infty}\Big|(\gamma+it)^{-b}\Big|\times \Big|\big(1-(\gamma+it)^{-1}\big)^{a-b}\Big|\d t.
\eeq

Note that the set $\{(\gamma+it)^{-1} :  t\in\mathbb{R} \} \subset \c$ is a circle with centre $(2\gamma)^{-1}$ and radius $(2\gamma)^{-1}$. Therefore the set $\{1-(\gamma+it)^{-1} :  t\in\mathbb{R} \} \subset \c$ is a circle with centre $1-(2\gamma)^{-1}$  and radius $(2\gamma)^{-1}$. Recall that $\gamma=1+\epsilon>0$, therefore this last circle does not touch the vertical line $\i\mathbb{R}$ and we have
 \beq\label{LastLemma_D_estimate}
 D=\max_{t\in\mathbb{R}}\Big\{|{\textnormal{arg}}(1-(\gamma+it)^{-1})|\Big\}<\frac{\pi}{2}.
\eeq
At the same time, we have for all $t\in \r$
\beqq
\frac{\epsilon}{1+\epsilon}\le \big|1-(\gamma+\i t)^{-1} \big|=\sqrt{\frac{\epsilon^2+t^2}{\gamma^2+t^2}}<1.
\eeqq
The above two estimates and the equality $|u^{v}|=|u|^{\re(v)}e^{|{\textnormal{arg}}(u) \times \im(v)|}$ (which is valid for all $u\in \c$ and $v\in \c$ with $\re(u)>0$, $\re(v)>0$) 
show that for all $t\in \r$, we have
\beq\label{LastLemma4}
\Big|\big(1-(\gamma+it)^{-1}\big)^{a-b} \Big | \le C_2(a_{0},a_{1},b,\epsilon)e^{D|Im(a)|},
\eeq
where 
\beqq
 C_2(a_{0},a_{1},b,\epsilon)=\max\left\{1,\left(\frac{\epsilon}{1+\epsilon}\right)^{a_0-\re(b)},\left(\frac{\epsilon}{1+\epsilon}\right)^{a_1-\re(b)}
\right\}.
\eeqq
Using (\ref{LastLemma3}) and (\ref{LastLemma4}), we conclude that
\beqq
\big|{}_1F_{1}(a,b,-z)\big|&\leq& C_1(b)C_2(a_{0},a_{1},b,\epsilon)e^{\epsilon z+ D|Im(a)|}\int_{-\infty}^{\infty} \big|(\gamma+it)^{-b}\big|\d t\\
&=&C(a_{0},a_{1},b,\epsilon)e^{\epsilon z+ D|Im(a)|}.
\eeqq
Note that the integral appearing in the above estimate converges since $\re(b)\in (1,2)$. This proves (\ref{ConluentAst}) for $z\ge 1$. 

Assume next that $z\in(0,1)$. Using \eqref{LastLemma1} with $\gamma=(1+\epsilon)/z$ and changing variables in the integral $w \mapsto (1+\epsilon+\i t)/z$ we get
\beqq
{}_1F_{1}(a,b,-z)=\frac{e^{1+\epsilon-z}}{2\pi}
\Gamma(b)\int_{\r}e^{\i t}
\left(1+\epsilon+\i t\right)^{-b}\left(1-\frac{z}{1+\epsilon+\i t}\right)^{a-b}dt.
\eeqq
Now we can proceed as in the case when $z>1$ noting that the set $\{1-z(1+\epsilon+\i t)^{-1} : t\in\mathbb{R}\}\subset \c$
is a circle with centre $1-z(2(1+\epsilon))^{-1}$ and radius $z(2(1+\epsilon))^{-1}$. 
As $0<z<1$ one can see that this only improves all the estimates above. 
For example, the estimate \eqref{LastLemma_D_estimate} also holds true and for all $t\in\r$, we have 
\beqq
\frac{\epsilon}{1+\epsilon}<\frac{1+\epsilon-z}{1+\epsilon}\le \big|1-z(1+\epsilon+\i t)^{-1} \big|=\sqrt{\frac{(1+\epsilon-z)^2+t^2}{(1+\epsilon)^2+t^2}}<1.
\eeqq
Therefore (\ref{ConluentAst}) is also true for $z\in(0,1)$.

Finally, we consider the case when $\re(b)\in (0,1)$. One can see that this case follows easily from the already established result valid 
for $\re(b) \in (1,2)$ and the following identity for the confluent 
hypergeometric function 
\beqq
b{}_1F_{1}(a,b,-z)=a{}_1F_{1}(a+1,b+1,-z)+(b-a){}_1F_{1}(a,b+1,-z),
\eeqq
see formula 9.212.3 in \cite{Jeffrey2007}.
\end{proof}

\noindent{\it Proof of Proposition \ref{AppProp2}.}
The equation (\ref{AppLemma2-2}) follows from  (\ref{eqn_M_parabolic_cylinder}) and the fact that 
\beqq
D_{-s}(0)=\frac{\sqrt{\pi}2^{-\frac{s}2}}{\Gamma\lb\frac12(s+1)\rb} ,
\eeqq
see the definition of the parabolic cylinder function (\ref{def_parabolic_cylinder}). 

Let us prove (ii). Assume first that $s\in \c$ is a fixed number which satisfies $\re(s)\in\lb\frac{1}{4},\frac{3}{4}\rb$. Using Lemma \ref{BmLemma1} and Fubini Theorem we find that the Mellin transform of $k(x)$ is given by
\beq\label{MF1}
\mathcal{M}_{\mu,1}(s)=\iint\limits_{\r_+^2}
F(s,y,z) 
 \p(J_{1}\in \d y, J_2\in \d z),
\eeq
where
\beqq
F(s,y,z)=\frac{1}{\sqrt{2\pi z}}\int_{0}^{\infty}x^{s-1}e^{-\frac{(x-\mu y)^{2}}{2z}}\d x. 
\eeqq
According to Lemma \ref{AppLemma1}, since $\re(s)\in\lb\frac{1}{4},\frac{3}{4}\rb$  the Mellin transform of $F(s,y,z)$ in the $z$-variable exists for all $w$ such that
$0<\re(w)<1/8$ and is given by
\beq\label{def_Gsyw}
G(s,y,w):=\int_0^{\infty} F(s,y,z) z^{w-1} \d z=2^{w-1}(-\mu y)^{2w+s-1}\frac{\Gamma(s)\Gamma(1-2w-s)}{\Gamma(1-w)}.
\eeq
Using Lemma \ref{Gamma_estimates}, we find that for every $s$ such 
that $1/4<\re(s)<3/4$ there exists $C=C(s)>0$ such that for all $w$ with $\re(w)=1/16$ we have
\beq\label{estimate_G}
\left|G(s,y,w)\right|< C |y|^{\re(s)-\frac{7}{8}}e^{-|\im(w)|}.
\eeq
Therefore as $\left|G(s,y,w)\right|$ is absolutely integrable along the line $w=\frac{1}{16}+\i\r$ then $F(s,y,z)$ can be written as an inverse Mellin transform
\beqq
F(s,y,z)=\frac{1}{2\pi \i} \int_{\frac1{16}+\i \r} G(s,y,w) z^{-w} \d w.
\eeqq
From the above identity and (\ref{MF1}) we find that
\beqq
\mathcal{M}_{\mu,1}(s)=\frac{1}{2\pi \i}\int_{y=0}^{\infty}\int_{z=0}^{\infty}\int_{w \in \frac1{16}+\i \r} G(s,y,w) z^{-w} \d w\p(J_{1}\in \d y, J_2\in \d z).
\eeqq
Due to (\ref{estimate_G}) and the fact that $\e[J_1^{\re(s)-7/8}J_2^{-1/16}]<\infty$ since $1/4<\re(s)<3/4$ (see Lemma \ref{ProbaLemma}), we conclude that the function
$G(s,y,w) z^{-w}$ is absolutely integrable with respect to the measure $\d w \times \p(J_{1}\in \d y, J_2\in \d z)$. Thus, we can apply  Fubini's Theorem to the right-hand side of the above equation and with the help of (\ref{def_Gsyw}), we obtain
\beq\label{eqn_M_before_change}
\mathcal{M}_{\mu,1}(s)=\frac{\Gamma(s)}{2\pi \i} \int_{\frac1{16}+\i \r} \frac{\Gamma(1-2w-s)}{\Gamma(1-w)} \e\left[J_1^{2w+s-1}J_2^{-w}\right] 
2^{w-1} (-\mu)^{2w+s-1} \d w.
\eeq

Next, we perform a change of variables $w \mapsto \frac12 ( 1-u-s)$ (recall that $s$ is a fixed number) and obtain from (\ref{eqn_M_before_change})
\beq\label{eqn_M_after_change}
\mathcal{M}_{\mu,1}(s)=\frac{2^{-\frac12(s+1)}\Gamma(s)}{2\pi \i} \int_{\frac7{8}-s+\i \r} \frac{\Gamma(u)}{\Gamma(\frac12(u+s+1))} 
 \e\left[J_1^{-u}J_2^{\frac12(u+s-1)}\right] 
 (2\mu^2)^{-\frac{u}2} \d u.
\eeq
For $s$ fixed, such that $1/4<\re(s)<3/4$, we know from (ii) Lemma \ref{ProbaLemma} that  $\e\left[J_1^{-u}J_2^{\frac12(u+s-1)}\right]$ is a bounded analytic function of $u$ everywhere in the strip
$-1<\re(u)<17/18-\re(s)$ and hence bounded on $\re(u)=7/8-\re(s)$. The ratio of Gamma functions $\Gamma(u)/\Gamma(\frac12(u+s+1))$ decays exponentially (and uniformly) as $\im(u)\to \infty$ in the strip 
$-1<\re(u)<17/18-\re(s)$, and it has a unique simple pole at $u=0$, coming from $\Gamma(u)$. Thus we can shift the contour of integration in (\ref{eqn_M_after_change})
$7/{8}-s+\i \r \mapsto -1/2+\i \r$ and taking into account the residue at $u=0$ we finally obtain
\beq\label{eqn_M_after_shift}
\mathcal{M}_{\mu,1}(s)&=&\frac{2^{-\frac12(s+1)}\Gamma(s)}{\Gamma\lb\frac12(s+1)\rb } 
\e\left[J_2^{\frac12(s-1)}\right]\\ \nonumber
&+&
\frac{2^{-\frac12(s+1)}}{2\pi \i} \int_{-\frac1{2}+\i \r} \frac{\Gamma(s)\Gamma(u)}{\Gamma\lb\frac12(u+s+1)\rb } 
 \e\left[J_1^{-u}J_2^{\frac12(u+s-1)}\right] 
 (2\mu^2)^{-\frac{u}2} \d u.
\eeq
According to Lemma \ref{ProbaLemma}, $\e\left[J_1^{-u}J_2^{\frac12(u+s-1)}\right]$ is a bounded analytic function for $\re(w)=-1/2$ and $-1+\epsilon<\re(s)<1+\theta-\epsilon$, for any $\epsilon>0$. Due to Lemma \ref{Gamma_estimates}, the ratio of  
Gamma functions $\Gamma(u)/\Gamma(\frac12(u+s+1))$ decays exponentially as $\im(u) \to \infty$, $\re(u)=-1/2$ and 
uniformly in $s$ if 
$-1+\epsilon<\re(s)<1+\theta-\epsilon$. Therefore, the right-hand side in (\ref{eqn_M_after_shift}) defines a meromorphic function in the strip $-1<\re(s)<1+\theta$, 
which has a unique simple pole at $s=0$ (which comes from $\Gamma(s)$), and we can apply analytic continuation and conclude that (\ref{eqn_M_after_shift}) is valid for all $s$ in the 
strip $-1<\re(s)<1+\theta$. This ends the proof of part (ii).

Finally, let us prove (iii). Assume first that $0<\re(s)<1$. We use formulae (\ref{def_parabolic_cylinder}) and (\ref{eqn_M_parabolic_cylinder})
to find that 
\beq\label{proof_iii_n1}
\mathcal{M}_{\mu,1}(s)+\mathcal{M}_{-\mu,1}(s)=\Gamma(s)
\frac{2^{-\frac12(s-1)}}{\Gamma(\frac12(s+1))}
\e \left[ J_2^{\frac12 (s-1)} e^{-\frac{\mu^2}2 V} {}_1F_1\left(\frac{s}2,\frac12,\frac{\mu^2}2 V \right) \right].
\eeq
From the above formula and the identity $e^{-z}{}_1F_1(a,b,z)={}_1F_1(b-a,b,-z)$ (see formula (7) on page 253 in \cite{Erdelyi}), we conclude that
 (\ref{AppLemma2-3}) holds true for $0<\re(s)<1$.  Now our goal is to check that formula \eqref{AppLemma2-3} can be extended into the wider strip
$-1<\re(s)<1+\theta$. 

Assume that $\delta>0$ is a small number and that $-1+\delta<\re(s)<1+\theta-\delta$. It is clear that we can find $p=p(\delta)>1$ such that
for all $s$ in the strip $-1+\delta<\re(s)<1+\theta-\delta$, we have
$(\re(s)-1)p \in (-2,\theta)$. Define $q=p/(p-1)$. 
 According to Lemma \ref{ProbaLemma}, we can find $\epsilon>0$
small enough such that $\e \left[ \exp(\epsilon q \frac{\mu^2}{2} V)\right]<\infty$. Using Lemma \ref{LastLemma}, we see that 
there exists $D=D(\delta)\in (0,\pi/2)$ and $C=C(\delta)>0$ such that for all $s$ in the strip $-1+\delta<\re(s)<1+\theta-\delta$, we have
\beqq
\bigg |{}_1F_1\left(\frac{1-s}2,\frac12,-\frac{\mu^2}2 V \right) \bigg|< C e^{\epsilon\frac{\mu^2}2 V + D |\im\lb\frac{s}2\rb |}.
\eeqq
Therefore, we can use H\"older inequality with $p$ and $q$ defined as above and estimate the expectation in the right-hand side of (\ref{AppLemma2-3}) as follows
\beq\label{1F1_estimate_exp_decay}
\Bigg | \e \left[ J_2^{\frac12 (s-1)} {}_1F_1\left(\frac{1-s}2,\frac12,-\frac{\mu^2}2 V \right) \right] \Bigg|
< C e^{ D |\im\lb \frac{s}2\rb |}\left( \e\left[ J_2^{\frac12 (\re(s)-1) p} \right] \right)^{\frac1p}
 \left( \e\left[ e^{\epsilon q \frac{\mu^2}{2} V}  \right] \right)^{\frac1q}<\infty,
\eeq
where in the last step we have used the fact that $\frac12(\re(s)-1)p \in (-1,\theta/2)$. This shows that 
the expectation in the right-hand side of (\ref{AppLemma2-3}) is well-defined for all $s$ such that $-1+\delta<\re(s)<1+\theta-\delta$, and since
$\delta>0$ is an arbitrary  small number, we can  extend the validity of this equation into the whole strip $-1<\re(s)<1+\theta$. 
\qed
\vspace{0.5cm}

The next theorem is our second main result in this section and it opens the way for the application of powerful complex-analytical tools.

\begin{theorem}\label{thm_asymptotic_decay}
Assume that $\xi$ satisfies condition \eqref{assumption_first_moments}. 
For any $\mu\in \mathbb{R}$ and any small number $\delta>0$, there exist constants $A=A(\mu,\sigma,\delta)>0$ and $B=B(\mu,\sigma,\delta)>0$ such that
\begin{equation}\label{EstimateMT}
\big|\mathcal{M}_{\mu,\sigma}(s)\big|\leq A  e^{-B|Im(s)|},
\end{equation}
for all $s\in\mathbb{C}$ such that $\re(s)\in(-1+\delta,1+\theta-\delta)$ and $|\im(s)|>1$. 
\end{theorem}
\begin{proof}
Note that $I_{\mu,\sigma}\stackrel{d}=\sigma I_{\frac{\mu}{\sigma},1}$, hence $\mathcal{M}_{\mu,\sigma}(s)=\sigma^{s-1}\mathcal{M}_{\frac{\mu}{\sigma},1}(s)$, therefore without  loss of generality  we can 
assume $\sigma=1$.

Since $-1/2$ is not a pole for $\Gamma(s)$, we use Lemma \ref{Gamma_estimates} and conclude that there exists $C_1>0$ such that
for all $y\in \r$ 
\begin{equation}\label{GammaAsympt2}
\Big|\Gamma(-\frac{1}{2}+iy)\Big|\leq C_1e^{-\frac{7\pi}{16}|y|}.
\end{equation}
At the same time, from Lemma \ref{Gamma_estimates} we find that for all $x\in (-1,1+\theta)$ 
and $y\in \r$ there exists $C_2>0$ such that  
\begin{equation}\label{GammaAsympt1}
\Big|\Gamma(x+iy)\Big|\geq C_{2}e^{-\frac{9\pi}{16}|y|}.
\end{equation}
First let us assume that $\mu=0$. Then \eqref{EstimateMT} follows immediately from Lemma \ref{Gamma_estimates} and \eqref{AppLemma2-2} since $\left|\e \left[J^{\frac{s-1}{2}}_{2}\right]\right|<C(\delta)$, for $\re(s)\in(-1+\delta,1+\theta-\delta)$. The latter is obvious from Proposition \ref{Proposition1} for $J_{2}$.\\

Next, assume that $\mu<0$. Thanks to \eqref{AppLemma2-1} and Lemma \ref{ProbaLemma}, we get that
\begin{equation}\label{IntermediateEstimate}
\left|\mathcal{M}_{\mu,1}(s)\right|\leq \tilde{C}(\mu,\delta,\varepsilon)\lb\frac{\left|\Gamma(s)\right|}{|\Gamma(\frac{1}{2}(s+1))|}+
\left|\Gamma(s)\right|\int_{\r}
\frac{|\Gamma\left(-\frac12+\i y \right)|}
{\left|\Gamma\left(\frac{1}{2}\left(\frac12+\i y +s\right)\right)\right|}\d y\rb.
\end{equation}
From Lemma \ref{Gamma_estimates}, we deduce that for $-1<\re(s)<1+\theta$ and $|\im(s)|>1$ there exists $C_3>0$ such that
\beqq
\frac{\left|\Gamma(s)\right|}{|\Gamma(\frac{1}{2}(s+1))|}\leq C_3 e^{-\frac{\pi}{6}|\im(s)|},
\eeqq
which shows that the first term in \eqref{IntermediateEstimate} is decaying exponentially as $\im(s) \to \infty$. 
Next, from Lemma \ref{Gamma_estimates}, we know that for  $-1<\re(s)<1+\theta$ and $|\im(s)|>1$ there exists $C_4>0$ such that
\beqq
|\Gamma(s)|<C_4 e^{-\frac{7 \pi}{16} |\im(s)|}.
\eeqq
Using this fact and estimates \eqref{GammaAsympt2} and \eqref{GammaAsympt1}, we see that for $|\im(s)|>1$
\beqq
 \left|\Gamma(s)\right|\int_{\r}
\frac{|\Gamma\left(-\frac12+\i y \right)|}
{\left|\Gamma\left(\frac{1}{2}\left(\frac12+\i y +s\right)\right)\right|}\d y\leq
\frac{C_{1}}{C_2}\left|\Gamma(s)\right|\int_{-\infty}^{\infty}e^{-\frac{7\pi}{16}|y|+\frac{9\pi}{32}|y+\im(s)|}\d y
\\ \le
C_4\frac{C_{1}}{C_2} e^{-\frac{7\pi}{16} |\im(s)|} \int_{-\infty}^{\infty}e^{-\frac{7\pi}{16}|y|+\frac{9\pi}{32}|y|+\frac{9\pi}{32}|\im(s)|}\d y
 \\ \leq 
C_4\frac{C_{1}}{C_2}  e^{-\frac{5\pi}{32} |\im(s)|} \int_{-\infty}^{\infty}e^{-\frac{5\pi}{32}|y|} \d y =
 C_{5}e^{-\frac{5\pi}{32}|\im(s)|}.
\eeqq
The above estimate shows that the second  term in \eqref{IntermediateEstimate} is decaying exponentially as $\im(s) \to \infty$, which ends the proof in
the case $\mu<0$. 

Finally, let us consider the case when $\mu>0$. In the proof of part (iii) of Proposition \ref{AppProp2} (see inequality \eqref{1F1_estimate_exp_decay}), 
we have established that for every $\delta>0$
there exist constants  $D=D(\delta)\in (0,\pi/2)$ and $C=C(\delta)>0$ such that for all $s$ in the strip $-1+\delta<\re(s)<1+\theta-\delta$ we have
\beqq
\Bigg | \e \left[ J_2^{\frac12 (s-1)} {}_1F_1\left(\frac{1-s}2,\frac12,-\frac{\mu^2}2 V \right) \right] \Bigg|
< C e^{ D |\im(\frac{s}2)|}.
\eeqq
Then from (\ref{AppLemma2-3}), we find that
for all $s$ in the strip $-1+\delta<\re(s)<1+\theta-\delta$ 
\beq\label{bound_n1}
\big|{\mathcal{M}}_{\mu,1}(s) \big| <  \Bigg| \frac{\Gamma(s)}{\Gamma(\frac12(s+1))} \Bigg|
C e^{ D |\im(\frac{s}2)|}+\big|{\mathcal{M}}_{-\mu,1}(s) \big|.
\eeq
Using Lemma \ref{Gamma_estimates} and the fact that $D<\pi/2$,  we conclude that there exist $C_1>0$ such that for all 
$s$ in the strip $-1+\delta<\re(s)<1+\theta-\delta$ and $|\im(s)|>1$, we have
\beqq
\Bigg| \frac{\Gamma(s)}{\Gamma(\frac12(s+1))} \Bigg|<C_1 e^{-\frac{1}{4}\left( \frac{\pi}2+D \right)  |\im(s)|}.
\eeqq
Therefore the first term in (\ref{bound_n1}) can be bounded by
\beqq
C C_1 e^{-\frac{1}{4}\left( \frac{\pi}2-D \right)  |\im(s)|},
\eeqq
and it decays exponentially as $\im(s)\to \infty$ since $\pi/2-D>0$.  This ends the proof in the case $\mu>0$,
since we have already established that the second term in  (\ref{bound_n1}) decays exponentially to zero.
\end{proof}

\begin{corollary}\label{InfDiff}
 Assume that $\xi$ satisfies condition \eqref{assumption_first_moments}. The function $k(x)$ is infinitely differentiable on $\r \setminus \{0\}$ and $k(x)\in C^{1}(\r)$.
\end{corollary}
\begin{proof}
The fact that $k(x) \in C^1 (\r )$  was already established in Lemma \ref{BmLemma1}. 
 Assume that $x\in (0,\infty)$. Applying Mellin transform inversion, we find that 
 \beq\label{kx_Mellin_transform_inversion}
 k_{\mu,\sigma}(x)=\frac{1}{2\pi\i}\int_{\frac{1}{2}+\i \r}x^{-s}\mathcal{M}_{\mu,\sigma}(s)\d s,
 \eeq
where the integral converges absolutely since \eqref{EstimateMT} guarantees exponential decay of $\mathcal{M}_{\mu,\sigma}(s)$ on the line 
$1/2 + \i \r$. This exponential decay also guarantees that for every $n\geq 0$ the functions
 \[\prod_{i=0}^{n-1}(s+i)\mathcal{M}_{\mu,\sigma}(s)\]
 are absolutely integrable along the line $1/2 + \i \r$, which shows by differentiation under the integral in \eqref{kx_Mellin_transform_inversion} that $k_{\mu,\sigma}(x)\in C^{\infty}(0,\infty)$. 
 Noting that $-I_{\mu,\sigma}\stackrel{d}=I_{-\mu,\sigma}$ we deduce that $k_{\mu,\sigma}(x)\in C^{\infty}(-\infty,0)$.
\end{proof}

\begin{corollary}\label{cor_asymptotic_decay}
Assume that $\xi$ satisfies condition \eqref{assumption_first_moments} and that $\hat{\rho}>0$, $\theta>0$.  For any $\mu \in \r$ and any small number $\delta>0$ the estimate \eqref{EstimateMT}
holds uniformly in the strip $\re(s)\in (-1-\hat \rho+\delta,1+\theta-\delta)$.
\end{corollary}
\begin{proof}
The statement about the exponential decay follows from Theorem \ref{thm_asymptotic_decay}, the functional equation \eqref{RecurBr} and the fact that $\psi(z)=O(z^2)$
uniformly in the strip $\re(s)\in (-1-\hat \rho+\delta,1+\theta-\delta)$. The latter fact follows from 
\eqref{SpecialFormExponents} (see also Proposition 2 in \cite{Be}).
\end{proof}

Theorem \ref{thm_asymptotic_decay} is very important for several reasons. First of all, as we have seen in Corollary \ref{InfDiff}, 
it implies smoothness of $k(x)$ on $\r\setminus\{0\}$. This should be compared with the case $\sigma=0$, where it is known that 
$k(x)$ may be non-smooth on $(0,\infty)$. For example, if $\xi$ has bounded variation and negative linear drift $\mu_{\xi}$, then $k(x)$ may be non-smooth at point $-1/\mu_{\xi}$, see Proposition 2.1 in \cite{CPY} and remark  2 in \cite{Ku11}. Secondly, as we will see in the next section, Theorem  \ref{thm_asymptotic_decay} together with Theorem \ref{BmTheorem1} will allow us to use simple techniques from Complex Analysis, such as shifting the contour of integration in the inverse Mellin transform, to prove rather strong results about the asymptotic behaviour of $k(x)$ as $x\to 0^+$ or $x\to +\infty$.


\section{Applications}\label{section_applications}


In this section, we present several applications of the results obtained in the previous section. 
We are still working under the same assumptions as in Section \ref{section_BM_exp_functionals}, i.e. we consider the exponential functional
$I_{\mu,\sigma}$ defined by \eqref{ExpFunBM} under the assumptions: $\e[|\xi_1|]<\infty$, $\e[\xi_1]<0$ and $\sigma=\sigma_\eta>0$. 

Our main tools are the meromorphic extension of $\mathcal{M}_{\mu,\sigma}$, Tauberian theorems and Mellin inversion with shifting of the contour of integration. We will also use the functional equation \eqref{RecurBr} and the estimate \eqref{EstimateMT} developed in Section \ref{section_BM_exp_functionals}. In Theorem \ref{Asym0} we derive some asymptotic results for $k(x)$ as $x\rightarrow 0$, while in Theorem \ref{thm_asym_infty} we discuss the behaviour of  $\p(I_{\mu,\sigma}>x)$ and $k(x)$ as $x\rightarrow\infty$, thus strengthening significantly
some of the results of Lindner and Maller \cite[Theorem 4.5]{LindnerMaller2005} in this 
special case when $\eta_{s}=\mu s+\sigma B_{s}$. We note that under further assumptions much stronger results are within reach for the asymptotic behaviour of $\p(I_{\mu,\sigma}>x)$ and $k(x)$ both as $x\rightarrow 0$ and $x\rightarrow\infty$. In order to illustrate the techniques, 
we choose to work with a rather simple but nevertheless very useful for applications class of processes $\xi$ which have hyper-exponential jumps 
(see \cite{Cai2009127,CaiKou2010,KUKYPA10}).
The same results can be easily generalized to more general class of L\'evy processes with jumps of rational transform (see \cite{Ku11}).

Finally we point out that $\p(I_{\mu,\sigma}>x)$ can be associated to ruin probability for certain actuarial models, see for example Theorem 4 in \cite{BS08}.


\subsection{General results about asymptotic behaviour of $k(x)$}


Our first theorem in this section deals with the asymptotic behaviour of $k(x)$ at zero. As usual, we define the ``floor'' (or ``integer part'') function as $\lfloor x \rfloor = \max\{n \in {\mathbb Z}\; : \; n \le x\}$.  We recall that $\hat \rho$ is defined by \eqref{def_rho_hat_rho_and_theta}: if $\hat \rho>0$ then the
L\'evy measure of $\xi$ has exponentially decaying negative tail with the rate of decay equal to $\hat \rho$.

\begin{theorem}\label{Asym0}
Assume that $\xi$ satisfies condition \eqref{assumption_first_moments} and that $\theta>0$, $\hat{\rho}>0$. Then  for every integer $m\ge 0$ and 
$\epsilon \in (0,1)$ such that $m+\epsilon<1+\hat \rho$ we have  
\begin{equation}\label{FstOrdAsymp1}
k_{\mu,\sigma}(x)=\sum_{n=0}^{m}\frac{b_{n}}{n!}x^{n}+O(x^{m+\epsilon}),\;\;
\textnormal{ as }\; x\to 0,
\end{equation}
where the coefficients $\{b_n\}_{n\ge 0}$ are defined recursively: $b_{-1}=0$, $b_{0}=k(0)$ and 
\begin{equation}\label{ReccurentRelation1}
b_{n+1}=\frac{2}{\sigma^{2}}\lb \mu b_{n}-\psi(-n)b_{n-1} \rb, \;\;\; \textnormal{ for } 0\le n<\hat{\rho}.
\end{equation}
In particular, $k(x)\in C^{1+\lfloor \hat\rho \rfloor}(\r)$, and if $\hat{\rho}=\infty$ then $k(x)\in C^{\infty}(\r)$. Moreover, as Remark \ref{SmoothnessBreakdown} shows $k''(0)$ may fail to exist.
\end{theorem}
\begin{proof}
Recall from Corollary \ref{thm_asymptotic_decay} that $\mathcal{M}_{\mu,\sigma}(s)$ is analytic in the strip $\re(s) \in (-1-\hat \rho,1+\theta)$ and has simple poles at all negative integers $-n$ such that $0\le n<1+\hat{\rho}$. Define
\beq\label{def_residue_an}
a_n=a_{n}(\mu,\sigma)=\textnormal{Res}(\mathcal{M}_{\mu,\sigma}(s) \; : \; s=-n), \;\;\; 0\le n<1+\hat{\rho}.
\eeq
Choose $c<1+\hat \rho$, such that $c \notin {\mathbb{N}}$.
We start from the Mellin transform inversion formula \eqref{kx_Mellin_transform_inversion}, use the fact that 
$\mathcal{M}_{\mu,\sigma}(s)$ decays exponentially as $\im(s) \to \infty$ (and uniformly in $\re(s)$) and shift the contour of integration 
$1/2+\i \r \mapsto  -c + \i \r$ while taking into account the residues at points $-n$ to obtain
\beq\label{K_asymptotic_proof_n1}
k_{\mu,\sigma}(x)=\frac{1}{2\pi \i}\int_{s=\frac{1}{2}+\i\r}x^{-s}\mathcal{M}_{\mu,\sigma}(s)\d s=\sum_{0 \le n <c}a_{n} x^{n}+\frac{1}{2\pi \i}\int_{-c+\i\r}x^{-s}\mathcal{M}_{\mu,\sigma}(s)\d s.
\eeq
The integral term in the right-hand side of the above equation can be estimated as follows
\beqq
\Bigg |\int_{-c+\i\r}x^{-s}\mathcal{M}_{\mu,\sigma}(s)\d s \Bigg |
\le x^c \int_{\r} \Big | \mathcal{M}_{\mu,\sigma}(-c + \i t) \Big | \d t,
\eeqq
therefore this term is $O(x^c)$ as $x\to 0^+$. 

Let us derive a recurrence relation for the coefficients $a_n$. First of all,
from Lemma \ref{Lemma_Mellin_transform_analytic_continuation} we find that 
$a_0=k(0)$ (this fact is also obvious from \eqref{K_asymptotic_proof_n1}). Next, from the definition \eqref{def_residue_an} and the fact that all the poles 
are simple, we find that 
\begin{equation}\label{relation111}
\mathcal{M}_{\mu,\sigma}(s)= \frac{a_{n}}{s+n}+O(1), \;\;\; s\to -n.
\end{equation}
Using formula \eqref{relation111} and the functional equation \eqref{RecurBr}, we find that as $s\to 0$ we have
\beqq
\frac{\psi(s)}{s} \mm_{\mu,\sigma}(s+1) + \mu \frac{a_0}{s} + O(1)+\frac{\sigma^2}2 (-1) \frac{a_1}{s} + O(1)=0.
\eeqq
Due to the fact that $\psi(s)/s\to \psi'(0)=\e[\xi_1]<\infty$, as $s\to 0$, and that $\mm_{\mu,\sigma}(s+1)\to 1$, as $s\to 0$, we conclude that $\mu a_0 - \sigma^2 a_1/2=0$. Following the same steps and considering the functional equation \eqref{RecurBr}, as $s\to -n$, we find that the coefficients $a_n$ 
satisfy the recurrence relation
\beqq
\frac{\psi(-n)}{-n}a_{n-1}+\mu a_{n}+\frac{\sigma^{2}}{2}(-n-1)a_{n+1}=0, \;\;\; n\ge 1.
\eeqq
Therefore, if we define $b_{n}=b_{n}(\mu,\sigma)=n! a_{n}(\mu,\sigma)$ then from the above equation, we obtain the recurrence relation \eqref{ReccurentRelation1}.

Combining all the above results we see that we have established (\ref{FstOrdAsymp1}), but only in the one-sided sense $x\to 0^+$.  
Using the fact that $I_{\mu,\sigma}\stackrel{d}{=}-I_{-\mu,\sigma}$ and repeating the above arguments, we obtain
\beqq
k_{\mu,\sigma}(-x)=k_{-\mu,\sigma}(x)=\sum_{n=0}^{m}\frac{b_{n}(-\mu,\sigma)}{n!}x^{n}+O(x^{m+\epsilon}),\;\;
\textnormal{ as }\; x\to 0^+.
\eeqq
Clearly, \eqref{FstOrdAsymp1} would be true if $b_{n}(-\mu,\sigma)=(-1)^n b_{n}(\mu,\sigma)$. This fact can be easily verified: 
using the recurrence relation \eqref{ReccurentRelation1}, we check that $b_{n}(\mu,\sigma)$ is a polynomial in $\mu$ of degree $n$, which is odd (even) 
if $n$ is an odd (even) number. This ends the proof of asymptotic formula (\ref{FstOrdAsymp1}). 

Finally, formula (\ref{FstOrdAsymp1}) and the fact that $k(x) \in C^{\infty}(\r \setminus \{0\})$  imply that 
$k(x) \in C^{1+\lfloor \hat\rho \rfloor}(\r)$, which ends the proof of Theorem \ref{Asym0}.
\end{proof}

\begin{remark}\label{ThmSmallAsymRem1}
To the best of our knowledge this is the first general result on the behaviour of $k_{\mu,\sigma}(x)$ as $x\rightarrow 0$ in the case $\sigma>0$. At the same time there are several recent results concerning such behaviour when $\sigma=0$, see \cite{Ku11,KUPA11,PRS}. 
\end{remark}

Note that if the process $\xi$ is spectrally positive, or more generally, if $\Pi_{\xi}(\d x)$ restricted to $(-\infty,0)$ has exponential moments of arbitrary order, then $k(x) \in C^{\infty} (\r)$. 

Our next result provides an extensive account of the asymptotic behaviour of $\p(I_{\mu,\sigma}>x)$ as $x\rightarrow+\infty$.\\

\begin{theorem}\label{thm_asym_infty}
Assume that $\xi$ satisfies condition \eqref{assumption_first_moments} and that $\theta>0$ and  $\xi$ has a non-lattice distribution. 
\begin{itemize}
 \item[(i)] 
If one of the following conditions is satisfied 
\begin{itemize}
 \item[(a)] $\theta<\rho$,
\item[(b)] $\theta=\rho$, $\psi(\theta)=0$ and $\e[ \xi_1^2 \exp(\theta \xi_1)]<\infty$, 
\end{itemize}
then 
\beq\label{FstOrdAsymp2}
 \p(I_{\mu,\sigma}>x)=C x^{-\theta}+o\left(x^{-\theta}\right), \;\;\; x\to +\infty, 
\eeq
where $C=-R(\theta)/\theta$ and $R(\theta)$ is defined in \eqref{def_C}. 
\item[(ii)] If $\theta=\rho$ and $\psi(\theta)<0$ then
\beq\label{FstOrdAsymp}
 \p(I_{\mu,\sigma}>x)=o\left(x^{-\theta}\right), \;\;\; x\to +\infty. 
\eeq
\end{itemize}
Moreover, if $\mu \le 0$ then the asymptotic expressions \eqref{FstOrdAsymp2} and \eqref{FstOrdAsymp} can be differentiated and leads to an asymptotic expression for $k(x)$. 
\end{theorem}
\begin{proof}
Let us prove part (i). For $x\ge 0$ we define  
\beqq
 F(x):=\int\limits_0^x \p(I_{\mu,\sigma}>y^\frac{1}{2\theta}) \d y. 
\eeqq 
Using integration by parts in the same way as we did above when dealing with the Mellin transform of the $F_2(k;v)$ term in the proof 
of Theorem  \ref{BmTheorem1} (see also equation \eqref{proof_F2kv}),  we find that for all $s$ in the strip $\re(s)\in (0,\frac12)$
\beq\label{eqn_hat_Fs}
\hat F(s):=\int_0^{\infty} x^{s-1} \d F(x)=s^{-1} {\mathcal{M}}_{\mu,\sigma}(1+2\theta s). 
\eeq
 Due to Corollary \ref{Corollary11}, the function $\hat F(s)-C/(1/2-s)$ is continuous in the strip $0< \re(s) \le 1/2$, therefore we can apply 
Wiener-Ikehara Theorem (see Theorem 7.3 in \cite{BaDi}) and conclude that as $x \to +\infty$ 
\beqq
F(x)=2 C \sqrt{x}+o\left(\sqrt{x}\right).
\eeqq
Using the above asymptotic expression, the fact that  $\p(I_{\mu,\sigma}>y^\frac{1}{2\theta})$ is a decreasing function of $y$ and applying the Monotone Density Theorem we obtain \eqref{FstOrdAsymp2}.

The proof of part (ii) is very similar: now we use Corollary \ref{Corollary11} to find that $\hat F(s)$ is continuous in the strip $0< \re(s) \le 1/2$,
therefore by applying Wiener-Ikehara Theorem we conclude that as $x \to +\infty$, we have $F(x)=o(\sqrt{x})$, and applying the Monotone Density Theorem we 
obtain \eqref{FstOrdAsymp}.

If $\mu<0$ then from Lemma \ref{BmLemma1} we know that $k(x)$ is a decreasing function on $(0,\infty)$, 
therefore we can apply the Monotone Density Theorem to \eqref{FstOrdAsymp2} or \eqref{FstOrdAsymp} and obtain the corresponding 
asymptotic expression for $k(x)$. 
\end{proof}

\begin{remark}\label{ThAsympRem1}
Note that, for all $t>0$, 
\[I_{\mu,\sigma}=\int_{0}^{t}e^{\xi_{s}}\d\eta_{s}+e^{\xi_{t}}I'_{\mu,\sigma},\]
where $I'_{\mu,\sigma}\stackrel{d}=I_{\mu,\sigma}$ and $I'_{\mu,\sigma}$ independent of $\lb e^{\xi_{t}}, \int_{0}^{t}e^{\xi_{s}}\d\eta_{s}\rb $. Recalling that in \cite[Prop. 4.1\,;\,Theorem 4.5]{LindnerMaller2005} the authors use $-\xi$ for $\xi$ in the definition of $I_{\mu,\sigma}$, we point out that the authors supplement the theory of random recurrent equations developed in \cite{Goldie91} to deduce general results for the behaviour of $\p\lb I(\xi,\eta)>x\rb$. For the case when $\eta_{t}=\mu t+\sigma B_{t}$ their result translates to 
\[\lim_{x\to\infty}x^{\theta}\p\lb I_{\mu,\sigma}>x\rb=C_{+}\geq 0; \lim_{x\to\infty}x^{\theta}\p\lb I_{\mu,\sigma}<-x\rb=C_{-}\geq 0;\,\lim_{x\to\infty}x^{\theta}\p\lb \left|I_{\mu,\sigma}\right|>x\rb=C_{+}+C_{-}>0\]
under the conditions (following our notation) that either $\rho>\theta\geq 1$ or $\theta<1<\rho$. Our assumptions are much weaker, see  (i) Theorem \ref{thm_asym_infty} and we also compute the constants $C_{\pm}$. 

Moreover, assuming $\mu\leq 0$, or otherwise working with $I_{-\mu,\sigma}\stackrel{d}=-I_{\mu,\sigma}$, one can show that $-R(\theta)>0$, thus $C_{+}>0$ and hence $C_{+}+C_{-}>0$. To prove that 
$-R(\theta)>0$,  we consider two cases: when $\theta\geq 1$ this follows directly from \eqref{def_C} and the fact
that $\psi'(\theta)>0$, and when
$\theta \in (0,1)$ this follows from \eqref{def_C} and the fact that $(\theta-1)\mathcal{M}_{\mu,\sigma}(\theta-1)>0$ (the latter is true due to \eqref{analytic_continuation_2} and the fact that $k(x)$ is decreasing, see Lemma \ref{BmLemma1}).  

Finally we note that, despite dealing with the asymptotics of $I(\xi,\eta)$ for general $\eta$, the methodology in \cite{LindnerMaller2005} cannot seemingly be improved to yield stronger results for the special case when $\eta$ is a Brownian motion with drift.
\end{remark}
\begin{remark}\label{ThAsympRem2}
The case when $\sigma=0$ has been completely dealt with in \cite{MZ,RI1,r2007}. We note that the technique applied there again relies on the random recurrence equations studied in \cite{Goldie91} and the authors are able to obtain results in part (i), condition (b) under the weaker 
assumption $\e\left[ \xi_1 \exp(\theta \xi_1)\right]<\infty$. A recent paper by V. Rivero \cite{r2009} addresses the case when the process $\xi$ has convolution equivalent L\'evy measure, the main tools are fluctuation theory of L\'evy processes and an explicit path-wise representation
of the exponential functional.  
\end{remark}


\subsection{Case study: processes with hyper-exponential jumps}\label{CaseStudy}


In this section, we will show how our methods can be extended to derive quite strong results about the density of the exponential functional, 
provided that we impose additional restrictions on the L\'evy process $\xi$. In particular,  we will need more information about the analytical structure of the Laplace exponent $\psi(z)$. 
Our purpose in this section is not to prove the most general results possible, 
but rather to present the ideas and give the flavour of the results which can be derived.

Let us consider a simple (but very useful) class of processes having hyper-exponential jumps (see \cite{Cai2009127, CaiKou2010, KUKYPA10}). In this 
case the L\'evy measure of a process $\xi$ is essentially a mixture of exponential distributions
\beq\label{def_Pi_hyp_exp}
\Pi_{\xi}(\d x)={\mathbf{1}}_{\{x>0\}}\sum\limits_{n=1}^N a_n e^{-\rho_n x}\d x+
{\mathbf{1}}_{\{x<0\}}\sum\limits_{n=1}^{\hat N} \hat a_n e^{\hat\rho_n x}\d x.
\eeq
where all the constants $a_n, \hat a_n, \rho_n, \hat \rho_n$ are strictly positive. 
Since $\lambda=\Pi_{\xi}(\r)<0$, the process $\xi$ can be represented as Brownian motion with drift plus a compound Poisson process
\beqq
\xi_t=\mu_{\xi}t+\sigma_{\xi} W_t + \sum\limits_{n=1}^{N(\lambda t)} Y_i,
\eeqq
where $N(t)$ is the standard Poisson process and $Y_i$ are i.i.d. random variables with distribution $\p(Y_i \in \d y)=\lambda^{-1} \Pi_{\xi}(\d y)$. Note that $\mu_{\xi}$ is the linear drift of the process $\xi$, it is easy to relate it to the constant $b_{\xi}=\e[\xi_1]$ by 
$b_{\xi}=\mu_{\xi}+\lambda \e[Y_1]$. 

Formula \eqref{def_Pi_hyp_exp} implies that the Laplace exponent of a hyper-exponential process is a rational function of the form
\beq\label{psi_partial_fractions}
\psi(z)=\frac{\sigma_{\xi}^2}2 z^2+\mu_{\xi} z + z \sum\limits_{n=1}^N  \frac{a_n}{\rho_n(\rho_n-z)} 
- z \sum\limits_{n=1}^{\hat N}  \frac{\hat a_n}{\hat \rho_n(\hat\rho_n+z)}.
\eeq
For hyper-exponential processes, it is  known (see \cite{Cai2009127, KUKYPA10}) that the equation $\psi(z)=0$ has only real simple solutions. Denote the
positive solutions as
$\{\zeta_m\}_{1\le m \le M}$, where we assume that they are arranged in increasing order. It is also known that 
$M=N+1$ if (i) $\sigma_{\xi}>0$ or (ii) $\sigma_{\xi}=0$ and $\mu_{\xi}>0$, and $M=N$ otherwise 
(see \cite{Ku11}). Note that in our previous notation \eqref{def_rho_hat_rho_and_theta}, we have $\theta=\zeta_1$, $\rho=\rho_1$ and
$\hat \rho=\hat \rho_1$.  

Using this information about zeros and poles of $\psi(z)$ and the functional equation \eqref{RecurBr} it is easy to see that $\mathcal{M}_{\mu,\sigma}(s)$ can be extended to a meromorphic 
function, with poles at the points 
\beqq
\{\zeta_m+n \; : \; m\ge 1, \; 1\le n \le N\}\cup\{-\hat \rho_n-m \; : \; m\ge 1, \; 1\le n \le \hat N\}\cup
\{-m \; : \; m\ge 0\}.
\eeqq
If we further assume that 
\beq\label{assumptions_hyper-exp}
 \begin{cases}
 & \zeta_i-\zeta_j \notin {\mathbb{Z}} \; \textnormal{ for all } \; 1\le i,j \le M \; \textnormal{ and } \; i \ne j, \\
 & \hat\rho_i \notin {\mathbb{N}} \; \textnormal{ and } \; \hat \rho_i-\hat \rho_j \notin {\mathbb{Z}} \; \textnormal{ for all } \; 1\le i,j \le \hat N
\; \textnormal{ and } \; i \ne j, 
\end{cases}
\eeq
then it is clear that all the poles of $\mathcal{M}_{\mu,\sigma}(s)$ are simple. Let us introduce the following notations
\beqq
c_{i,j}&=&-\frac{1}{(\zeta_i)_j}{\textnormal{Res}}\lb\mathcal{M}_{\mu,\sigma}(s) \; : \; s=j+\zeta_i \rb , \;\;\; 1\le i \le M, \;\;\; j \ge 1,  \\
b_{i,j}&=&(1+\hat \rho_i)_j {\textnormal{Res}}\lb\mathcal{M}_{\mu,\sigma}(s) \; : \; s=-j-\hat \rho_i \rb , \;\;\; 1 \le i \le \hat N, \;\;\; j \ge 1,
\eeqq
recall that $(a)_n=a(a+1)\dots(a+n-1)$ denotes the Pochhammer symbol.
Our next goal is to compute coefficients $c_{i,j}$ and $b_{i,j}$ in terms of the Mellin transform $\mm_{\mu,\sigma}(s)$. 
Let us fix $i$ such that $1\le i \le M$. Since $\zeta_i$ is a simple root of a rational function $\psi(z)$, we have 
$\psi(z)=\psi'(\zeta_i)(z-\zeta_i)+O((z-\zeta_i)^2)$ as $z\to \zeta_i$. This fact and the functional equation \eqref{RecurBr}
show that
\beq\label{def_ci1}
c_{i,1}=\frac{1}{\psi'(\zeta_i)} \left( 
 \mu\mathcal{M}_{\mu,\sigma}(\zeta_i)+\frac{\sigma^{2}}2(\zeta_i-1)\mathcal{M}_{\mu,\sigma}(\zeta_i-1) \right).
\eeq
Next, using the functional equation \eqref{RecurBr} and the same technique as in the proof of Theorem \ref{Asym0}, we obtain a recursion equation 
\beq\label{recurrence_cij}
c_{i,j+1}=-
\frac{1}{\psi(j+\zeta_i)} \left( 
 \mu c_{i,j}+\frac{\sigma^{2}}2 c_{i,j-1} \right), \;\;\; j\ge 1,
\eeq
where we have defined $c_{i,0}=0$. Next, let us fix $i$ such that $1\le i \le \hat N$. Formula 
 \eqref{psi_partial_fractions} implies that $\psi(z)$ has a simple pole at $z=-\hat\rho_i$ with residue $\hat a_i$. Again, we use this fact and the functional 
equation \eqref{RecurBr} to conclude that 
\beq\label{def_bi1}
b_{i,1}=-\frac{2}{\sigma^2} \frac{\hat a_i}{\hat \rho_i} \mathcal{M}_{\mu,\sigma}(1-\hat \rho_i),
\eeq
and 
\beq\label{recurrence_bij}
b_{i,j+1}=\frac{2}{\sigma^2}
\left( \mu b_{i,j} -  \psi(-j-\hat \rho_i) b_{i,j-1} \right),  \;\;\; j\ge 1,
\eeq
where we have defined $b_{i,0}=0$. Recall that the coefficients
 $b_j=j! {\textnormal{Res}}(\mathcal{M}_{\mu,\sigma}(s) \; : \; s=-j)$ can be computed via the recurrence relation \eqref{ReccurentRelation1}

Our main result in this section is the following Theorem, which provides a complete asymptotic expansion for $k(x)$ 
as  $x\to 0^+$ and at $x\to +\infty$. The corresponding expansions as $x\to 0^-$ and $x\to -\infty$ can be obtained
by considering $k_{-\mu,\sigma}(x)=k_{\mu,\sigma}(-x)$.

\begin{theorem}\label{AsymptExp}
Assume that $\xi$ is a process with hyper-exponential jumps \eqref{def_Pi_hyp_exp}, which satisfies $\e[\xi_1]<0$ and for which conditions
\eqref{assumptions_hyper-exp} are satisfied. 
Then for every $c>0$ 
\begin{equation}
k_{\mu,\sigma}(x)=
\begin{cases}\label{eqn_asymptotics_hyperexp}
&\displaystyle
\sum_{0\le j < c}\frac{b_{j}}{j!}x^{j}+
\sum\limits_{i=1}^{\hat N} \sum\limits_{j\ge 1} 
{\mathbf{1}}_{\{ j+\hat \rho_i<c\}}  b_{i,j} \frac{ x^{j+\hat \rho_i}}{(1+\hat \rho_i)_j}  + O\left(x^{c}\right),
\;\;\; {\textnormal{ as }} \;\;  x \to 0^+, 
\\ \\
&\displaystyle
\sum\limits_{i=1}^M \sum\limits_{j\ge 1} {\mathbf{1}}_{\{ j+\zeta_i<c\}} c_{i,j} \frac{(\zeta_i)_j }{x^{j+\zeta_i}} + O\left(x^{-c}\right), 
\qquad \qquad \qquad\;\;\;\;\; {\textnormal{ as }} \;\; x\to +\infty.
\end{cases}
\end{equation}
\end{theorem}
\begin{proof}
From \eqref{psi_partial_fractions} it is clear that $\psi(z)=O(z^2)$ and $1/\psi(z)=O(1)$ as $\im(z)\to \infty$, $|\im(z)|>1$, and that these 
estimates are uniform in $\re(s)$. Therefore, using Theorem \ref{thm_asymptotic_decay} and the functional equation 
\eqref{RecurBr} we see that $\mathcal{M}_{\mu,\sigma}(s)$ decays exponentially as $\im(s) \to \infty$, $\im(s)>1$, and uniformly if $\re(s)$ belongs to a compact subset of $\r$. 
This shows that we can apply the same technique as in the proof of Theorem \ref{Asym0}: shift the contour of integration, collect all the residues 
and estimate the resulting integral. The details are left to the reader. 
\end{proof}

\begin{remark}\label{SmoothnessBreakdown}
 Note that if $\hat \rho=\hat \rho_1 \in (0,1)$ then the coefficient $b_{i,1}$ defined by \eqref{def_bi1} is strictly negative. Theorem \ref{AsymptExp}
shows that as $x\to 0^+$ we have 
\beqq
k(x) = k(0)+ k'(0) x + \frac{b_{i,1} }{1+\hat \rho} x^{1+\hat \rho}+o(x^{1+\hat \rho}),
\eeqq
 which implies that in this case
$k''(0)$ does not exist. 
\end{remark}

We would also like to point out that if conditions \eqref{assumptions_hyper-exp} are not satisfied, then $\mm_{\mu,\sigma}(s)$ will have multiple poles. This is
not a big problem, but it implies that the asymptotic expansions \eqref{eqn_asymptotics_hyperexp} will contain terms of the form $x^{\alpha} \ln(x)^k$,
where $-\alpha$ is the pole of $\mm_{\mu,\sigma}(s)$ and $k$ is a non-negative integer which is not greater than the multiplicity of the pole $-\alpha$. Also, results similar to Theorem 
\ref{AsymptExp} can be derived for a more general class of L\'evy processes, for example for processes which have jumps of rational transform, see 
\cite{Ku11} for results in the case $\sigma=0$.




\end{document}